\documentclass[twoside]{article}

\usepackage[utf8]{inputenc} 
\usepackage[T1]{fontenc}    
\usepackage{hyperref}       
\usepackage{url}            
\usepackage{booktabs}       
\usepackage{amsfonts}       
\usepackage{nicefrac}       
\usepackage{microtype}      
\usepackage{xcolor}         
\usepackage{mathtools}
\usepackage{amsthm}
\usepackage{amsmath}
\usepackage{amssymb}
\usepackage{graphicx}
\usepackage{dblfloatfix}

\usepackage[accepted]{aistats2024}
\usepackage[round]{natbib}

\DeclareMathOperator*{\vv}{\mathrm{vec}}
\global\long\def\eqdef{\overset{\text{def}}{=}}%
\global\long\def\E{\mathcal{E}}%

\global\long\def\P{\mathbf{P}}%
\global\long\def\J{\mathbf{J}}%

\global\long\def\R{\mathbb{R}}%

\newtheorem{theorem}{Theorem}[section]

\newtheorem{lemma}[theorem]{Lemma}

\theoremstyle{definition}
\newtheorem{defn}{Definition}[section]

\usepackage{authblk}

\DeclarePairedDelimiterX{\inp}[2]{\langle}{\rangle}{#1, #2}

\begin{document}

%

%

\twocolumn[

\aistatstitle{Fast and Accurate Estimation of Low-Rank Matrices from Noisy Measurements via Preconditioned Non-Convex Gradient Descent}

\aistatsauthor{ Gavin Zhang \And Hong-Ming Chiu \And  Richard Y. Zhang }

\aistatsaddress{Department of Electrical and Computer Engineering, University of Illinois-Urbana Champaign} ]

\begin{abstract}

Non-convex gradient descent is a common approach for estimating a low-rank $n\times n$ ground truth matrix from noisy measurements, because it has per-iteration costs as low as $O(n)$ time, and is in theory capable of converging to a minimax optimal estimate. However, the practitioner is often constrained to just tens to hundreds of iterations, and the slow and/or inconsistent convergence of non-convex gradient descent can prevent a high-quality estimate from being obtained. Recently, the technique of \emph{preconditioning} was shown to be highly effective at accelerating the local convergence of non-convex gradient descent when the measurements are noiseless. In this paper, we describe how preconditioning should be done for noisy measurements to accelerate local convergence to minimax optimality. For the symmetric matrix sensing problem, our proposed preconditioned method is guaranteed to locally converge to minimax error at a linear rate that is immune to ill-conditioning and/or over-parameterization. Using our proposed preconditioned method, we perform a 60 megapixel medical image denoising task, and observe significantly reduced noise levels compared to previous approaches.

\end{abstract}

\begin{figure}[t!]
	\centering
	\includegraphics[width = \columnwidth]{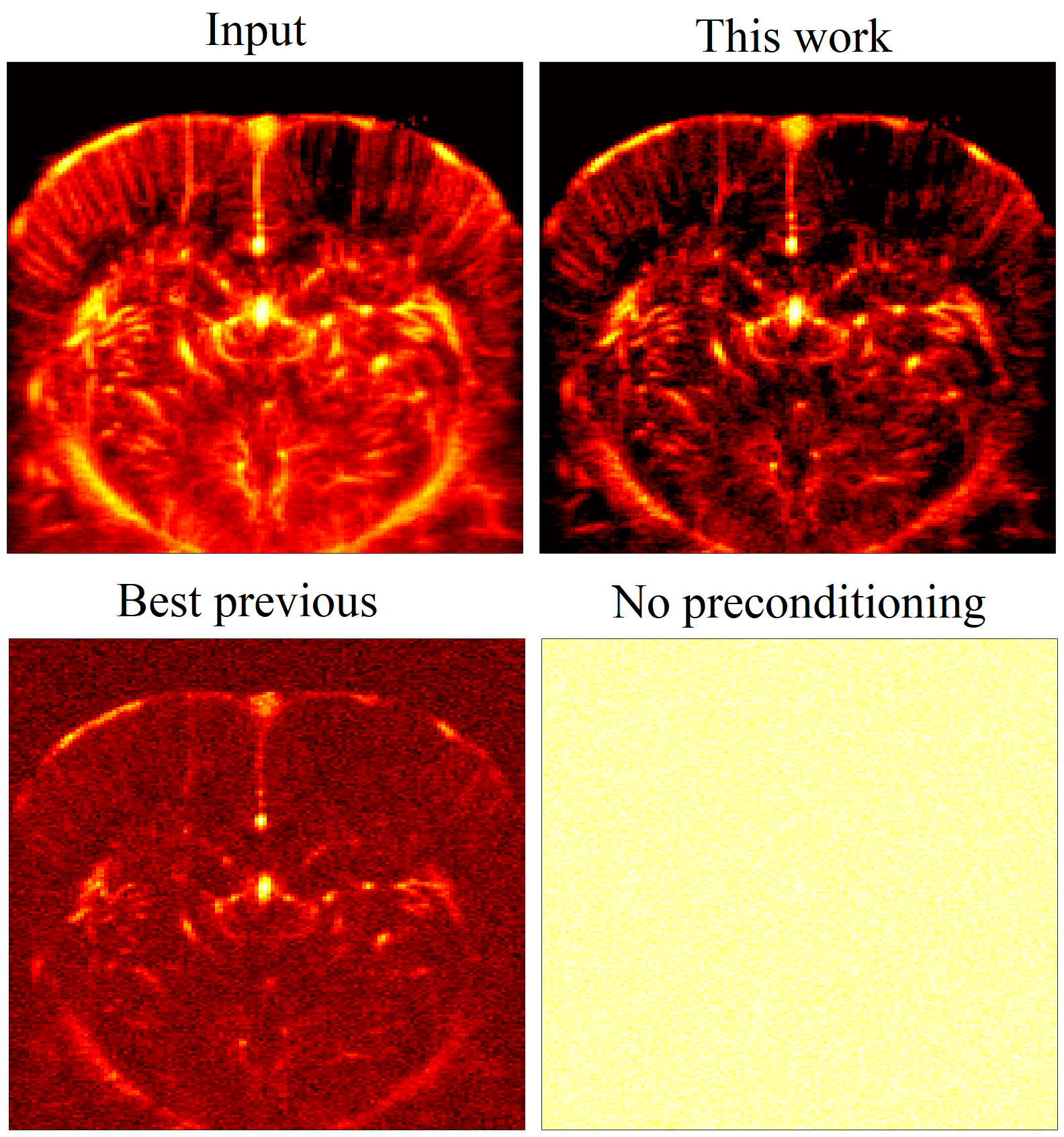}
	\caption{\small\textbf{Preconditioned gradient descent for a 60 megapixel medical image denoising task.} We denoise a 2400-frame ultrafast ultrasound image of a rat rain ($200\times 130$ pixels per frame) by running 30 iterations of the low-rank denoising procedure in \cite{demene2015spatiotemporal}. 
 Top-left: original noisy input. Top-right: image denoised and reconstructed by our preconditioning scheme in \eqref{alg:1}. Bottom-left: image obtained via the preconditioning scheme in \cite{zhang2021preconditioned}, which is the previous state-of-the-art. Bottom-right: image obtained by naive non-convex gradient descent without preconditioning.} 
    \label{fig:ultrasound}
\end{figure}

\section{INTRODUCTION}


We consider the \emph{low-rank matrix recovery} problem, which seeks to recover an $n\times n$ ground truth matrix $M^\star$ of low rank $r^\star$, given a small number $m$ of measurement matrices $A_i$ and noisy observations 
$y_i = \inp{A_i}{M^\star} + \varepsilon_i$, for indices $i\in\{1,\dots,m\}$. The main challenge lies in the fact that the presence of measurement noise $\varepsilon_i$ reduces the ``amount of useful information'' contained within the observations; it usually becomes impossible to recover $M^\star$ exactly. 
Instead, one aims to compute a \emph{minimax optimal} estimate $M\approx M^\star$, which is roughly defined as
the closest estimate of the ground truth $M^\star$ for the worst-possible scenario \cite{candes2010matrix,candes2011tight}. Informally, a minimax optimal estimate is the best achievable given the finite amount of useful information contained with the noisy observations \cite{lehmann2006theory}.


Minimax optimal estimations are highly desirable for real-world applications of low-rank matrix recovery. In medical imaging, for example, minimax optimality would assure the highest possible level of reconstruction accuracy, in order to minimize the chances of diagnostic errors, detect subtle changes or anomalies, and reduce the need for repeated scans or reanalysis. 
In response, an extensive body of literature has been developed over the past two decades on techniques for solving low-rank matrix recovery to guaranteed minimax optimality; see our detailed literature review in Section~\ref{subsec:related} below. 
Unfortunately, despite significant progress, existing state-of-the-art algorithms still have trouble achieving minimax optimality on many real-world applications, due to two perennial difficulties.

The first perennial difficulty is the enormous scale of real-world datasets. 
Today, the most common approach is \emph{non-convex gradient descent} (\cite{zheng2015convergent,zhao2015nonconvex,tu2016low,sun2016guaranteed}), or factored gradient descent (\cite{chen2015fast,park2017non,park2018finding}), which is to factor a candidate estimate $M=UV^T,$ and to directly optimize over its $n \times r$ low-rank factor matrices $U,V$, as in
\begin{equation}
\label{eq:obj}
\min_{U,V\in\mathbb{R}^{n\times r}} f(U,V) = \frac1m \sum_{i=1}^m (y_i - \inp{A_i}{UV^T})^2
\end{equation}
using an iterative local optimization algorithm like gradient descent
\begin{align}
\begin{split}\label{eq:gd}
U_{\text{new}} &= U - \alpha \nabla_U f(U,V), \\
V_{\text{new}} &= V - \alpha \nabla_V f(U,V),
\end{split}
\end{align}
in which $\alpha \in (0,1]$ is the step-size / learning rate. With a small rank $r\ll n$, each iteration
costs as low as $O(m+n)$ time and memory, and so in principle, the approach can scale to arbitrarily large values of $m$ and $n$. But real-world datasets routinely have $m$ and $n$ on the order of tens to hundreds of millions, and in practice, even a single iteration can take many minutes to several hours. 


The considerable expense of performing even a single iteration often constrains the practitioner to just a few tens to low hundreds of iterations. But this further exacerbates the second perennial difficulty, which is the inconsistent and sometimes extremely slow convergence of non-convex gradient descent. While the method is known to converge to minimax optimality given a \emph{sufficiently large} number of iterations (\cite{chen2015fast}), for many real-world datasets it is unable to do so with a \emph{reasonable} number of iterations. 

\subsection{Accelerating convergence via preconditioning}

Recently, there has been exciting progress on the use of \emph{preconditioning} to accelerate the local convergence of non-convex gradient descent for low-rank matrix recovery (\cite{mishra2012riemannian,tong2020accelerating,zhang2021preconditioned,zhang2022accelerating,xu2023power,zhang2023preconditioned}). This line of work is motivated by the observation that real-world ground truth matrices $M^\star$ often have excessively large condition numbers $\kappa=\lambda_1(M^\star)/\lambda_r(M^\star)$. In particular, if the search rank is over-parameterized as $r>r^\star$ with respect to the unknown true rank $r^\star$, then the condition number even diverges as $\kappa\to\infty$. Non-convex gradient descent is known to locally converge with a linear convergence rate like $\rho = 1- c/\kappa$ with absolute constant $c>0$ (\cite{zheng2015convergent,tu2016low}), and therefore experiences a significant slow-down with an excessively large or even diverging $\kappa$. In both cases, suitable preconditioning was shown to restore the linear convergence rate back to $\rho = 1- c$ (\cite{tong2020accelerating,zhang2021preconditioned}), as if the condition number were perfect $\kappa=1$.

However, the existing literature on preconditioning 
primarily focuses on the \emph{noiseless} instance of low-rank matrix recovery, which assumes $\varepsilon_i =0$ for all $i$. Indeed, it remains unclear how preconditioning should be done in the presence of measurement noise $\varepsilon_i \ne 0$. Experimentally, existing preconditioning methods for the noiseless case do not consistently accelerate convergence in the presence of noise; the acceleration is either lost at a much coarser error level than minimax optimal, or the iterates sporadically diverge. (See our experiments in Section~\ref{sec:experiment}) 

The existence of a significant gap between the noiseless and noisy cases is less surprising if we consider the underlying mechanism that allow preconditioning to work in the first place. Intuitively, preconditioning works by inverting one ill-conditioned matrix against another ill-conditioned matrix, in order to ``cancel out'' their ill-conditioning and obtain a well-conditioned matrix. This mechanism necessarily requires a precise alignment between the two matrices; in the presence of noise, slight ``misalignments'' can nullify the cancellation and render the preconditioning ineffective, or at worst even amplify the noise and cause divergence to occur.


\subsection{Our contribution: How to precondition in the presence of measurement noise?}

Our primary goal in this paper is to provide a \emph{principled} way to perform preconditioning on non-convex gradient descent, that is both effective and reliable for real-world applications of \emph{noisy} low-rank matrix recovery. To this end, we propose the following iterations for solving \eqref{eq:obj}:
\begin{align}
\begin{split}\label{alg:1}
U_{\text{new}} &= U - \alpha \nabla_U f(U,V) (V^T V + \eta I)^{-1}, \\
V_{\text{new}} &= V - \alpha \nabla_V f(U,V) (U^T U + \eta I)^{-1}, \\
\eta_{\text{new}} &= \beta \eta,
\end{split}
\end{align}
where $\alpha \in (0,1]$ is the step-size / learning rate as before in \eqref{eq:gd}, and $\beta \in [0,1)$ is a geometric decay rate for the regularization parameter $\eta$. Here, each gradient $\nabla_U f(U,V)$ and $\nabla_V f(U,V)$ is an $n\times r$ matrix, so the equation \eqref{eq:gd} says that each $r\times r$ matrix preconditioner $(V^T V + \eta I)^{-1}$ and $(U^T U + \eta I)^{-1}$ should be applied as a right matrix-matrix product onto the respective gradient. It is easily verified that the additional overhead of computing and applying the preconditioner is $O(r^3+nr^2)=O(n)$ time and $O(r^2)=O(1)$ memory, again assuming a small rank $r\ll n$.

The basic form of the preconditioned iterations \eqref{alg:1} is reminescent of previous work on noiseless low-rank matrix recovery (~\cite{mishra2012riemannian,mishra2016scaled,tong2020accelerating,zhang2021preconditioned,zhang2022accelerating,zhang2023preconditioned,xu2023power}). In the noisy setting, we provide strong theoretical and empirical evidence to argue that the most principled way to adjust the regularization parameter $\eta$ is to make it decay geometrically. This is in contrast to prior work, that either set $\eta=\sqrt{f(U,V)}$ at each iteration~(\cite{zhang2021preconditioned}), or simply fix $\eta$ to a constant for all iterations (\cite{xu2023power}). Our main message in this paper is that choosing $\eta$ correctly has a significant and outsized impact on the quality of acceleration in the noisy setting; the previous choices deliver an acceleration only at coarser error levels, and could even cause divergence. As shown in Figure~\ref{fig:ultrasound}, our method significantly improves upon the previous state-of-the-art on a real-world instance of low-rank matrix recovery arising in medical image denoising.

Rigorously, we prove, for the \emph{symmetric matrix sensing} instance of low-rank matrix recovery, that the preconditioned iterations in \eqref{alg:1} with a geometrically decaying $\eta$ locally converges to minimax optimality, at an accelerated linear rate that is immune to ill-conditioning and over-parameterization. It was previously shown that,
if the measurements $A_1,\dots,A_m$ satisfy the restricted isometry property (RIP) and that the measurement noise come from a zero-mean Gaussian $\varepsilon_i \sim \mathcal{N}(0,\sigma^2),$ then non-convex gradient descent with rank $r=O(r^\star)$ converges to an estimate $M$ with Frobenius norm error $\|M-M^\star\|_F=O(\sigma^2 n  r^\star \log n)$ (\cite{chen2015fast,zhuo2021computational}), which is indeed minimax optimal up to log factors (\cite{candes2011tight}). However, the actual convergence rate can be dramatically slowed by ill-conditioning and over-parameterization, to be as slow as sublinear (\cite{zhuo2021computational,zhang2021preconditioned}). A variant of \eqref{alg:1} known as PrecGD was proposed in \cite{zhang2021preconditioned} to accelerate convergence to minimax error, but this prior work required either perfect knowledge of the noise variance $\sigma^2$ (which is clearly unreasonable in practice), or a complicated and expensive cross-validation procedure to estimate the noise variance $\sigma^2$.

Our main result is that, under the same setting as the above, the preconditioned iterations in \eqref{alg:1} with geometrically decaying $\eta$ is guaranteed to converge to minimax optimal error, at the same local convergence rate as non-convex gradient descent with a perfect condition number $\kappa=1$. In particular: (i) the accelerated convergence rate is maintained all the way down to the minimax optimal error level of $O(\sigma^2 n r^\star)$; (ii) the acceleration is applicable to all initial points within a neighborhood of the ground truth. 
To the best of our knowledge, our result is the first to rigorously guarantee both two properties. 
Our analysis reveals a simple mechanism that explains the inconsistent performance of previous preconditioners in the noisy setting: the accelerated convergence is only maintained until the current error norm  $\|UV^T-M^\star\|_F$ reaches the same order of magnitude as the regularization parameter $\eta$, but an excessively small $\eta$ can actually cause the iterations to diverge. Therefore, the most natural and principled way to set $\eta$ is to allow it to geometrically decay alongside the error norm $\|UV^T-M^\star\|_F$, which is exactly what we proposed in \eqref{alg:1}.

\subsection{Limitations}

There are two main limitations with our theoretical analysis. First, we make two idealized assumptions (via the symmetric matrix sensing problem) that may not be satisfied in real-world datasets: (i) the underlying ground truth $M^\star$ is symmetric positive definite; (ii) that the measurements $A_i$ satisfy RIP. 
Here, we emphasize that the purpose of our theoretical analysis is to provide a \emph{rigorous justification} for the geometric decay of the regularization parameter $\eta$; in particular, it is not to guarantee performance on real-world datasets, as these will rarely satisfy the idealized assumptions (like RIP and incoherence) needed for a theoretical analysis to be possible. The symmetry assumption is primarily made to simplify our presentation; it can be mechanically overcome by repeating the analyses in \cite{park2017non,tong2020accelerating}, but we expect our conclusions to transfer largely verbatim to the non-symmetric case. We leave the extension of our theoretical results to the non-symmetric RIP setting as future work, and emphasize that our large-scale experiments on real-world datasets are indeed performed for a non-symmetric dataset that does not satisfy RIP.

Second, our theoretical analysis focuses on local convergence given a sufficient good initialization. For the sake of an end-to-end guarantee, we initialize our theoretical analysis using the standard technique of spectral initialization (\cite{tu2016low, chen2021spectral}). In practice, the enormous scale of real-world datasets often constrains the practitioner to just a few tens to low hundreds of iterations, so \emph{heuristic warm starts} are widely used to maximize the effectiveness of these few iterations (\cite{bercoff2011ultrafast, zhang2019spurious}). It is important to point out that slow convergence remains a critical issue even when a high-quality heuristic warm start is provided, as further progress towards minimax optimality is slowed by the slow convergence of the iterative algorithm. In this regard, our work in this paper answers the practical question: ``given a heuristic warm start, how do we refine this warm start to the best accuracy possible, while using as few iterations as possible?''

\subsection{\label{subsec:related}Related Work}

\paragraph{Non-convex gradient descent converges to minimax optimality}

For a wide range of problems relating to low-rank matrix recovery, if we are given a warm-start solution,
local refinements via GD is often capable of converging towards the ground truth. This has been shown rigorously for matrix sensing (\cite{tu2016low, zheng2015convergent, charisopoulos2021low}), matrix completion (\cite{sun2015nonconvex, jain2013low, chen2020nonconvex}), phase retrieval (\cite{candes2015phase, netrapalli2013phase, ma2018implicit}) and other related problems (\cite{li2019rapid, yi2016fast, chen2021bridging}).

For matrix sensing in particular, both  \cite{chen2015fast} and \cite{zhuo2021computational} showed that gradient descent achieves a statistical error of $O(\sigma^2 nr \log n )$, where $r$ is the search rank. Under the assumption that $r=O(r^\star)$, this error matches the minimax error noted in \cite{candes2011tight} up to log factors. In this work we prove that our method converges to the same statistical error as both \cite{chen2015fast, zhuo2021computational}, under exactly the same assumptions. In other words, our method does preconditioning without amplifying the statistical error at all. Under the warm-start setting, the question of whether the assumption $r=O(r^\star)$ is necessary for achieving minimax error is an open question even without preconditioning, and we do not attempt to resolve it in this work.



\paragraph{Accelerating local convergence via preconditioning}

The basic idea to precondition the gradient against $(V^T V)^{-1}$ and $(U^T U)^{-1}$ was first suggested in \cite{mishra2012riemannian}, and its convergence properties for the noiseless case were later studied in detail in \cite{tong2020accelerating} resulting in a method known as ScaledGD. Its extension to SGD was first proposed in \cite{mishra2016scaled}, and studied in \cite{zhang2022accelerating}. The idea to regularize with an identity perturbation and precondition against $(V^T V + \eta I)^{-1}$ and $(U^T U + \eta I)^{-1}$ was first suggested in \cite{zhang2021preconditioned} as a means to counteract the effects of over-parameterization $r>r^\star$, which resulted in a method known as PrecGD. A similar regularization was subsequently studied in \cite{xu2023power}.

All of these methods, as well as our proposed method, are able to overcome the slow convergence of non-convex gradient descent in the the noiseless setting. However, previous methods do not provide theoretical guarantees in the noisy setting. Empirically, their behaviors are inconsistent under noise. Specifically,  the ScaledGD of \cite{tong2020accelerating} can diverge in the case $r>r^*$. While this divergence is avoided by adding a regularization parameter as in \cite{zhang2013restricted} and \cite{xu2023power}, the regularization parameter itself can cause these methods to stagnate at a higher noise level, as seen in our experimental section (Figures \ref{fig:gaussian_matrix_sensing1} and \ref{fig:gaussian_matrix_sensing2}). In contrast, our method is the only one that maintains its acceleration all the way down to minimax optimality.

\paragraph{Small random initialization}
While most early work in non-convex gradient descent focused on local refinement of a warm-start initialization, a separate line of recent work focused on using a \textit{small} random initialization (\cite{li2018algorithmic, stoger2021small, ma2021implicit, ding2021rank, jin2023understanding}). For matrix sensing in particular, global convergence of GD was first  proven in \cite{li2018algorithmic} in the case $r=n$, and later refined in  \cite{stoger2021small} for the general over-parameterized case. Similar results have also been obtained in the asymmetric case (\cite{jiang2023algorithmic, chou2023induce}). For preconditioned methods, a similar analysis has been done in \cite{xu2023power} for a variant of the ScaledGD \cite{tong2020accelerating}. In the noisy setting, \cite{ding2022validation} first showed that GD with small random initialization converges to the minimax error by extending the analysis in \cite{stoger2021small}. A major strength of their theoretical analysis is that it no longer requires the assumption $r = O(r^*)$.

However, we emphasize that these theoretical results for small random initialization are not directly comparable to the results in this work. First, we provide theoretical guarantees for all initializations close to the ground truth. In contrast, small initialization relies on tracing a very specific and rapidly converging trajectory. In fact, we believe that this is the main reason that small random initialization achieves minimax optimality without requiring $r = O(r^*)$. However, in order to trace this specific trajectory, small random initialization forces an already good initial solution to be thrown away. In our experience, this means that GD has to use many more iterations to get back the warm-start that was thrown away (see Figure \ref{fig:gaussian_matrix_sensing3}).

\paragraph{Notations}
We use $\|\cdot\|_F$ to denote the Frobenius norm and $\|\cdot\|$ to denote the spectral norm of a matrix. We use $\inp{A}{B} = \mathrm{tr}(A^TB)$ to denote the standard matrix inner product. 
We use $\lesssim$ to denote an inequality that hides a constant factor. For a scalar function $f:\mathbb{R}^{n\times r}\to\mathbb{R}$, the gradient $\nabla f(X)$ is a matrix of size $n\times r$. For any matrix $M$, the eigenvalues and singular values are denoted by $\lambda_i(M)$ and $\sigma_i(M)$, arranged in decreasing order.

\section{MAIN RESULTS}

In our theoretical analysis, we consider the variant of low-rank matrix recovery known as \emph{symmetric matrix
sensing}, which aims to recover a positive semidefinite, rank-$r^{\star}$
ground truth matrix $M^{\star}\succeq 0$, from a small number $m$ of possibly
noisy measurements $y=\mathcal{A}(M^{\star})+\varepsilon$, where the linear measurement operator $\mathcal{A}$ is defined
\[
\mathcal{A}(M^{\star})=[\langle A_{1},M^{\star}\rangle,\langle A_{2},M^{\star}\rangle,\dots,\langle A_{m},M^{\star}\rangle]^{T}.
\]
Without loss of generality, we assume that the measurement matrices $A_i$ are symmetric. In addition, we will adopt the standard assumption that the unknown measurement noise modeled via the
length-$m$ vector $\epsilon$ is normally distributed
\[
\epsilon_i \sim \mathcal{N}(0,\sigma^2) \text{ for all } i\in \{1,\dots,m\}
\]
and that $\mathcal{A}$ satisfies the \textit{restricted isometry property }(RIP) \cite{candes2008restricted}. 

\begin{defn}[Restricted Isometry]
The linear operator $\mathcal{A}$ satisfies RIP with parameters
$(2r,\delta)$ if there exists constant $0\le\delta<1$ 
such that, for every rank-$2r$ matrix $M$, we have 
\begin{equation}
\label{def:rip}
(1-\delta)\|M\|_{F}^{2}\leq\|\mathcal{A}(M)\|^{2}\leq(1+\delta)\|M\|_{F}^{2}.
\end{equation}
\end{defn}

Specifically, we will always assume throughout the paper that $\mathcal{A}$ satisfies
RIP with parameters $(2r,\delta)$. We note that the RIP assumption is in line with existing work on the statistical optimality of gradient descent (\cite{chen2015fast,zhuo2021computational}) and preconditioned gradient descent (\cite{tong2020accelerating,zhang2021preconditioned}).  Under the mild assumption $r = O(r^*)$, it is also in line with prior work on convex methods (\cite{candes2008restricted,candes2010matrix,candes2011robust}) and small random initialization (\cite{li2018algorithmic, stoger2021small, ma2021implicit, ding2021rank, jin2023understanding, xu2023power}).

Given a warm-start close to the ground truth, our goal is to refine this warm-start by using gradient descent. In particular, we want to minimize the non-convex loss function in  \eqref{eq:obj} up to minimax optimal error. 

Since the ground truth $M^\star$ and the measurement matrices $A_i$'s are symmetric, if both the left and right factors $U, V$ in \eqref{alg:1} start at the same initial point, they will always stay the same. Therefore, if we denote $X = U= V$, with $X \in \mathbb{R}^{n\times r}$, then the iterations simplify to 
\begin{align}
\begin{split}\label{alg:new}
X_{\text{new}} &= X - \alpha \nabla f(X) (X^T X + \eta I)^{-1} \\
\eta_{\text{new}} &= \beta \eta
\end{split}.
\end{align}

Now we are ready to state our main result, which says that our algorithm always converges linearly to the minimax optimal error, at a linear rate that is affected by neither ill-conditioning nor over-parameterization. 

In particular, in equation \eqref{alg:new}, let the initial regularization $\eta_0$ satisfy  $\eta_0 \geq 2\sqrt{f_c(X_0)}$, and let the decay rate $\beta$ satisfy  $1>\beta \geq \sqrt{1-\frac{\mu}{4L}}$. Here $f_c(X_0)$ is the noiseless function value defined in \eqref{eq:fc}, and $\mu, L$ are both constants depending only on $r, r^*$ and $\delta$, which we define rigorously in Appendix A.2. Then we have the following result. 

\begin{theorem}
\label{thm:main}
Suppose that the initial point $X_{0}$ satisfies
$\|\mathcal{A}(X_{0}X_{0}^{T}-M^{*})\|^{2}<\rho^{2}(1-\delta)\lambda_{r^{*}}(M^{\star})^{2}$ with a radius $\rho>0$ that satisfies $\rho^{2}/(1-\rho^{2})\le(1-\delta^{2})/2$. Let the step-size $\alpha$ satisfy 
$\alpha \leq 1/L,$ where $L>0$ is a constant that only depends on $\delta$. At the $t$-th iteration, with high probability, we have 
\begin{align*}
\|X_tX_t^T-M^\star\|_F^2 \lesssim \max\left\{\beta^{2t}\cdot \|X_0X_0^T-M^\star\|^2_F, \mathcal{E}_{opt} \right\},
\end{align*} 
where $\mathcal{E}_{opt} = \frac{\sigma^{2}nr\log n}{m}$. Here the inequality $\lesssim$ hides a constant that only depends on $\delta$. 
\end{theorem}
A complete proof of Theorem \ref{thm:main} is presented in the appendix. In the next section, we sketch out the key ideas behind its proof. First, we make a few important observations. 

In Theorem \ref{thm:main}, we require an initial point that satisfies $\|\mathcal{A}(X_{0}X_{0}^{T}-M^{*})\|^{2}<\rho^{2}(1-\delta)\lambda_{r^{*}}(M^{\star})^{2}$. This requirement is standard and appeared in all previous works on preconditioned methods (\cite{tong2020accelerating, zhang2021preconditioned, zhang2023preconditioned}). The only exception is \cite{xu2023power}, which uses a small random initialization. In practice, a common way to obtain such an initial point is through domain specific heuristics. In ultrafast ultrasound (\cite{bercoff2011ultrafast}) for instance, the noisy version of the ultrasound image itself can serve as a warm-start, since it is already close to the ground truth. However, even without heuristics, such an initial point can be achieved using spectral initialization (see Proposition 6 of \cite{zhang2021preconditioned}), in which we simply need to compute one SVD factorization. 


In addition, we also need a decay rate that satisfies $1>\beta \geq \sqrt{1-\frac{\mu}{4L}}$. Although $\mu$ and $L$ are in general hard to estimate, we find that in practice $\beta$ is extremely robust. In our experiments, any value of $\beta$ satisfying $0.5\leq \beta<1$ was sufficient for linear convergence. 

We also note that in Theorem \ref{thm:main}, the convergence rate is crucially independent of the condition number $\kappa$, since $\mu$ and $L$ has no dependence on $\kappa$. In addition,  the statistical error that our algorithm converges to exactly matches that of \cite{chen2015fast} and \cite{zhuo2021computational}, which proved that GD converges to an error of $O(\sigma^2 nr \log n)$. In other words, our preconditioner exponentially accelerates GD without amplifying the noise at all. Under the mild assumption $r = O(r^*)$, this rate matches the minimax rate in \cite{candes2011tight}.


\section{KEY IDEA and PROOF SKETCH}
\label{sec:keyidea}


The recent work of \cite{zhang2021preconditioned} proposed a preconditioned variant of gradient descent called PrecGD to restore its linear convergence rate:
\begin{gather*}
X_{t+1} = X_t-\alpha \nabla f(X_t) (X_t^TX_t + \eta_tI)^{-1},\\    
C_{lb} \|X_tX_t^T-M^\star\| \leq \eta_t \leq  C_{ub} \|X_tX_t^T-M^\star\|,
\end{gather*}
where $C_{lb},C_{ub}$ are fixed constants.
To understand why our algorithm achieves minimax optimality and immunity to ill-conditioning and over-parameterization all at the same time, it is instructive to first see how PrecGD can fail in the noisy case.

To maintain linear convergence, the key contribution of \cite{zhang2021preconditioned} is the crucial observation that the regularization parameter $\eta_t$ must be within a constant factor of the error $\|X_tX_t^T-M^\star\|$. In the noiseless case, simply setting $\eta_t = \sqrt{f(X_t)}$ will imply $\eta_t = \Theta(\|X_tX_t^T-M^\star\|)$. However, in the noisy setting finding the right $\eta_t$ requires an accurate estimate of the noise variance, which is in general very difficult. 


\subsection{Key Innovations}

Maintaining the right amount of regularization of is the most important ingredient for our method to succeed. The regularization $\eta_t$ used in PrecGD has to be perfect because linear convergence 
requires two contradictory properties to intersect, namely gradient dominance (also known as the PL-inequality) and Lipschitz smoothness. When $\eta_t$ is too large, gradient dominance is lost. When $\eta_t$ is too small, Lipschitz smoothness is lost. The analysis in \cite{zhang2021preconditioned} suggests that the choice of $\eta_t$ is extremely fragile and delicate.

Surprisingly, we find that this is not the case. In fact, we will show that the choice of $\eta_t$ is not delicate, but rather robust. Our method avoids the need to choose the optimal regularization parameter altogether by simply letting $\eta$ decay with some rate $\beta<1$. It turns out that this extremely simple choice of the regularization parameter will \textit{automatically} maintain the right amount of regularization needed for linear convergence due to a phenomenon we call ``coupling''. 

This phenomenon can be intuitively understood as a race in which the two runners $\eta_t$ and $\mathcal{E}_t {=} \|X_tX_t^T-M^\star\|_F$, are connected using a rubber band. When $\eta_t$ and $\mathcal{E}_t$ begin to grow apart, the rubber band will exert a counteracting force and pull them back together. As a result, the amount of regularization is always right. This happens because that $\eta_t$ itself controls how fast $\mathcal{E}_t$ decays.  If the regularization parameter $\eta_t$ is large compared to $\mathcal{E}_t$, then our algorithm behaves more like gradient descent. As a result, our algorithm briefly stagnates, allowing $\eta_t$ to catch up and become close to $\mathcal{E}_t$ again. Similarly, if $\eta_t$ is small, our algorithm begins to converge faster. Thus, the error $\mathcal{E}_t $ decays quickly, and will eventually catch up to $\eta_t$. 

This coupling of the regularization parameter and the error is precisely why we can avoid the expensive procedure used in PrecGD to estimate the noise variance and approximate $\mathcal{E}_t$. We \textit{implicitly} maintain the right amount of regularization, so that our algorithm always converges linearly, even in ill-conditioned, over-parameterized, and noisy settings.

\subsection{Proof Sketch}
In this section we sketch the main steps of the proof of our main result, Theorem \ref{thm:main}, and defer the full proof to the appendix. Our proof consists of two components: the first is the observation that the PL-inequality, which is lost in the case $r>r^\star$, can be restored under a change of norm, as long as the preconditioner $P = (X_t^TX_t + \eta_t \cdot I)^{-1}$ has the 	``correct'' amount of regularization $\eta_t$. The second component is the observation that $\eta_t$ and $\mathcal{E}_t =\|X_tX_t^T-M^\star\|$ are coupled together, meaning that they can never be too far apart. 

\begin{figure*}[b!]
	\centering
    \vspace{-1em}
	\includegraphics[width = 0.5\textwidth]{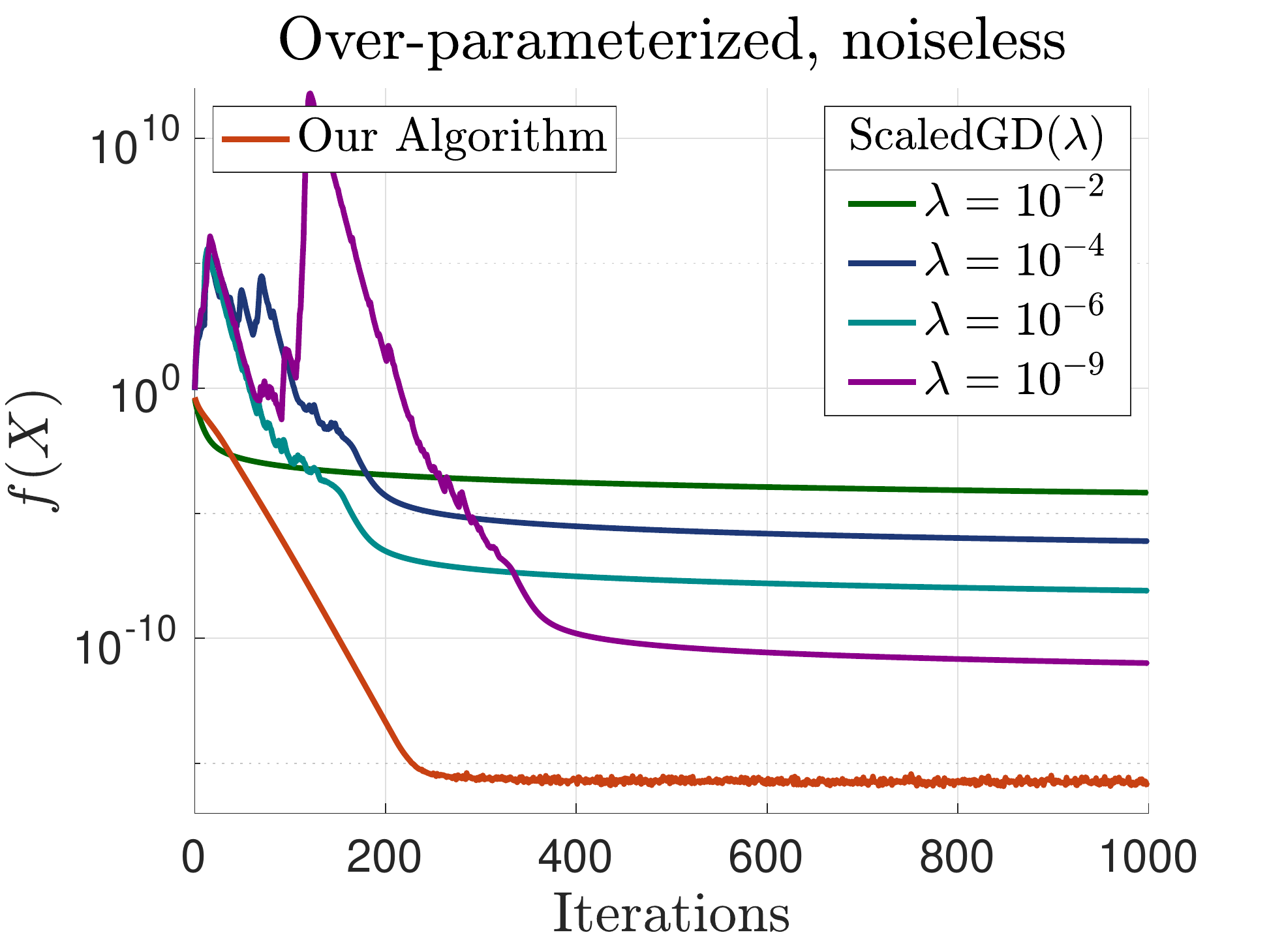}%
	\includegraphics[width = 0.5\textwidth]{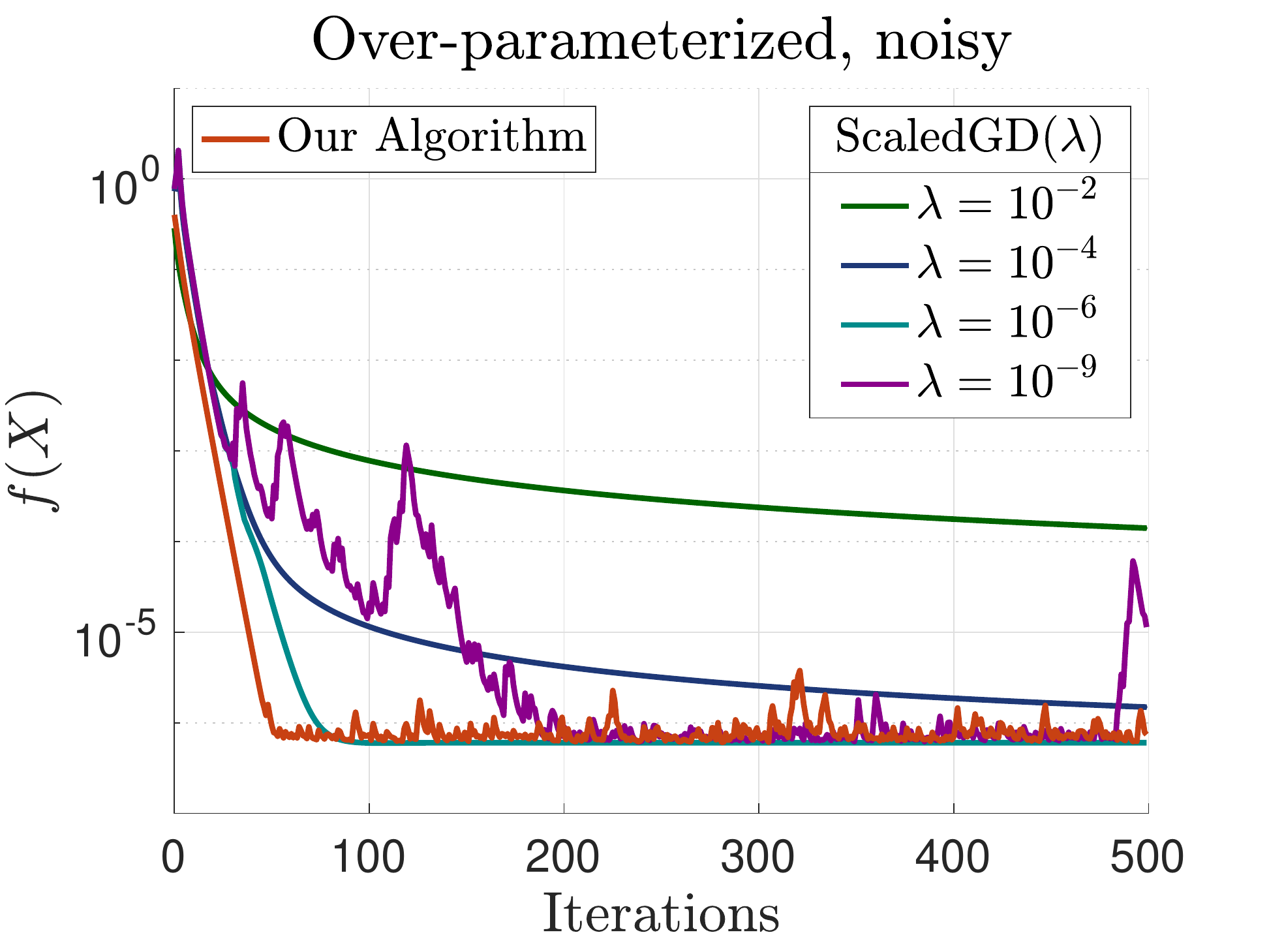}
    \vspace{-2em}
	\caption{\textbf{Convergence of our algorithm and ScaledGD$(\lambda)$ using spectral initialization}. Left: Noiseless measurements. Right: Noisy measurements with noise variance $\sigma=10^{-6}$.}
    \label{fig:gaussian_matrix_sensing1}
\end{figure*}
\begin{figure*}[t!]
	\centering
    \vspace{-0.5em}
	\includegraphics[width = 0.5\textwidth]{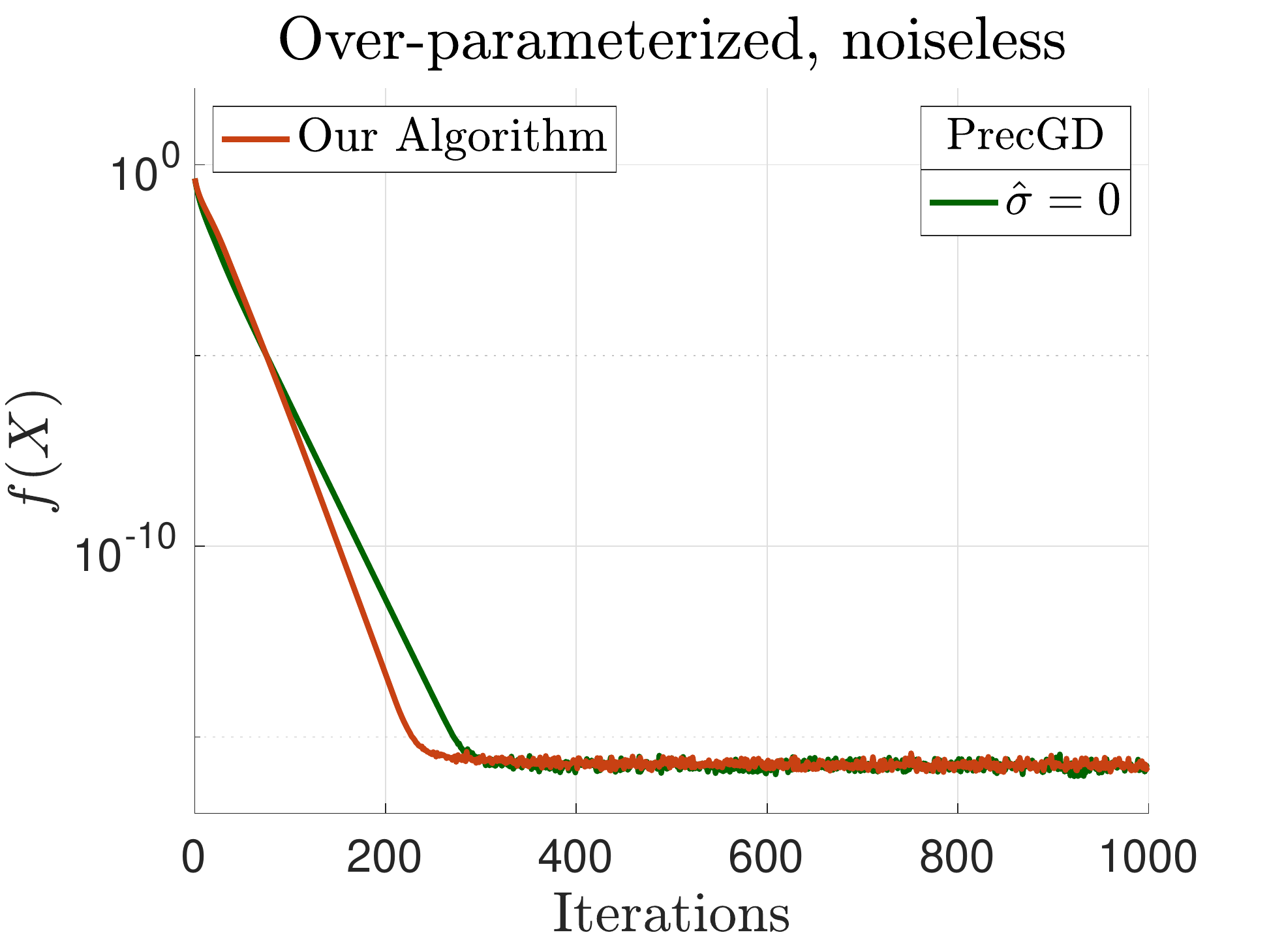}%
	\includegraphics[width = 0.5\textwidth]{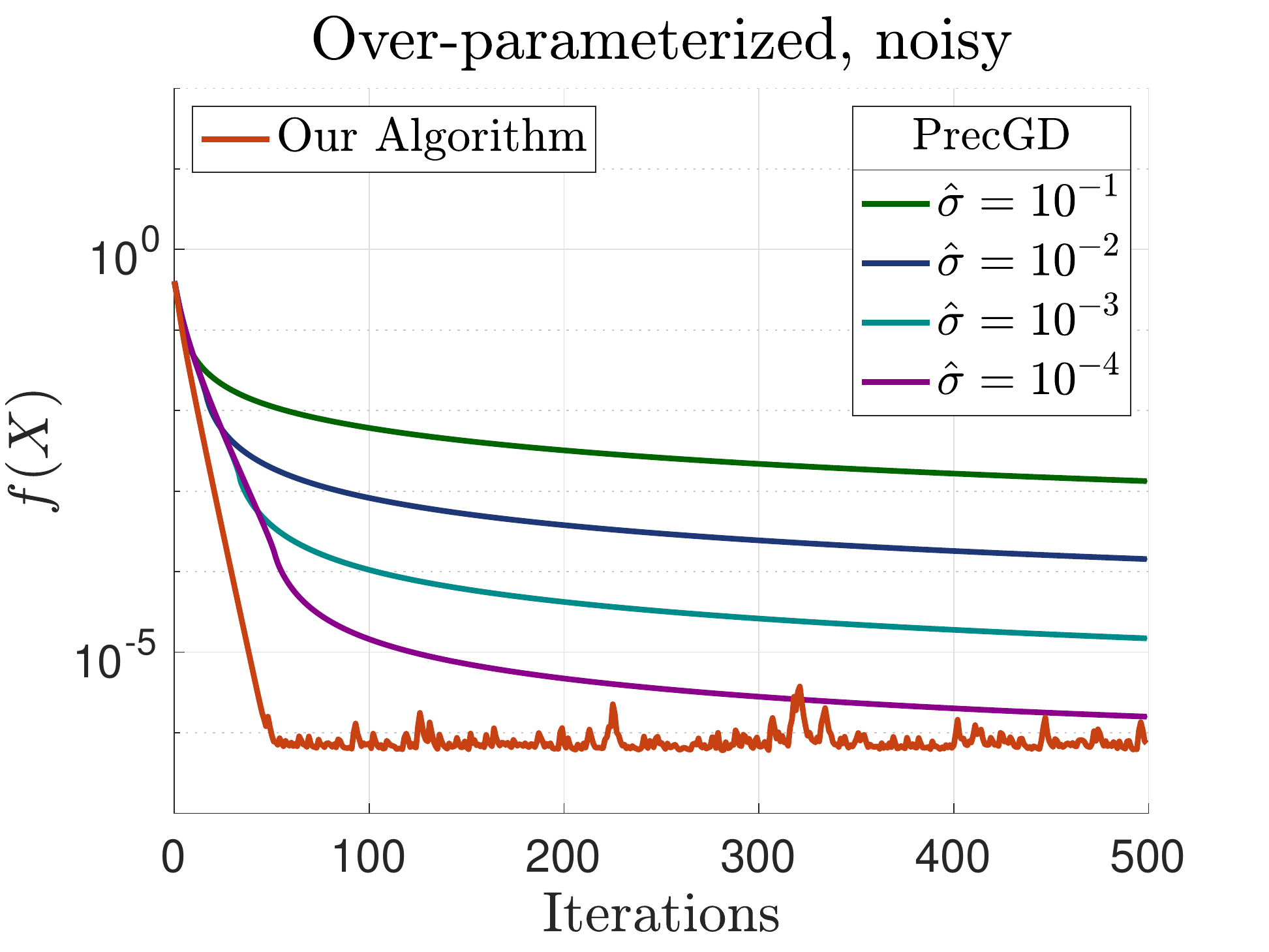}
    \vspace{-2em}
	\caption{\textbf{Convergence of our algorithm and PrecGD using spectral initialization}. Left: Noiseless measurements. Right: Noisy measurements with noise variance $\sigma=10^{-6}$.}
    \label{fig:gaussian_matrix_sensing2}
\end{figure*}
\begin{figure*}[t!]
	\centering
	\includegraphics[width = 0.5\textwidth]{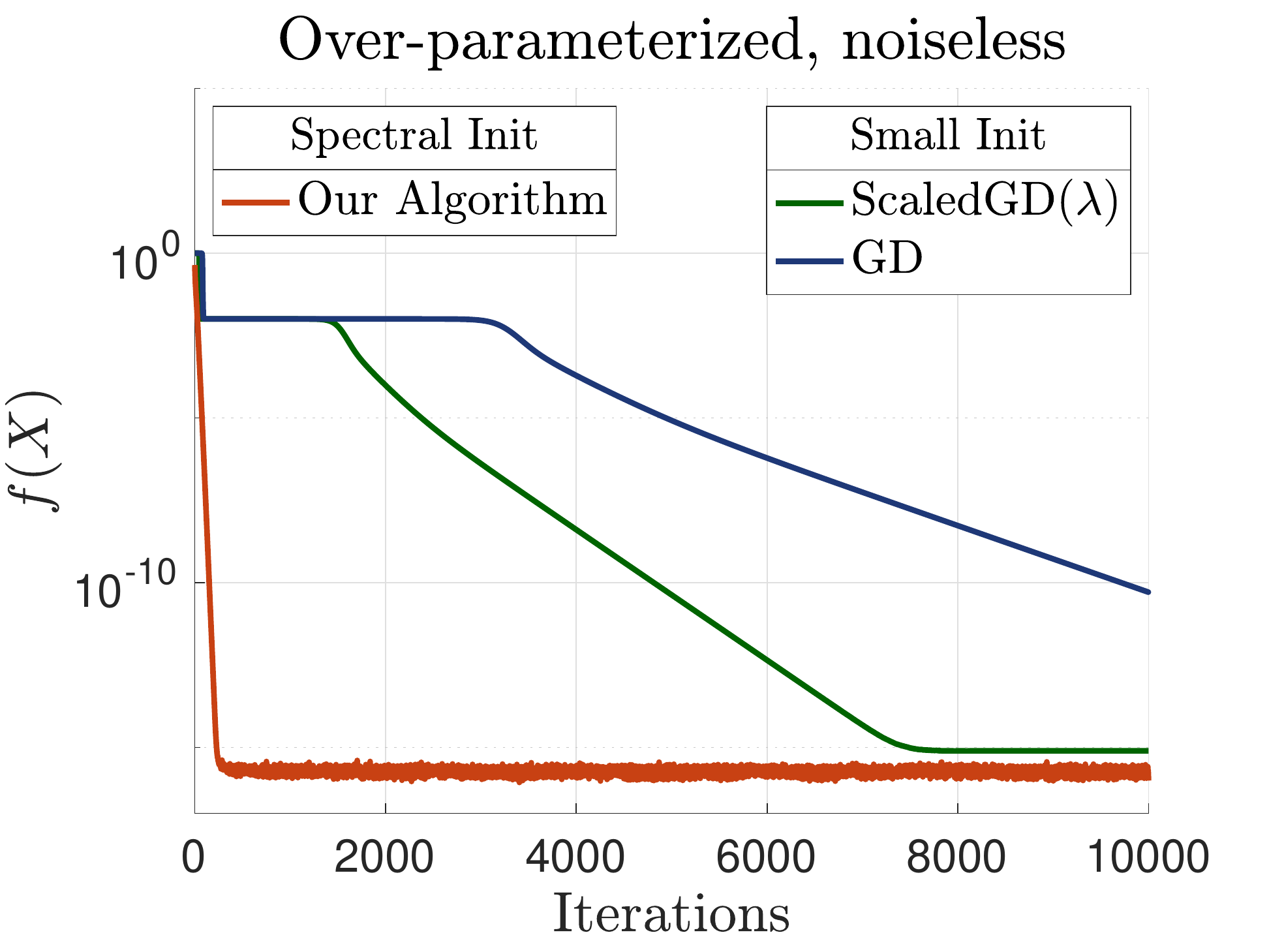}%
	\includegraphics[width = 0.5\textwidth]{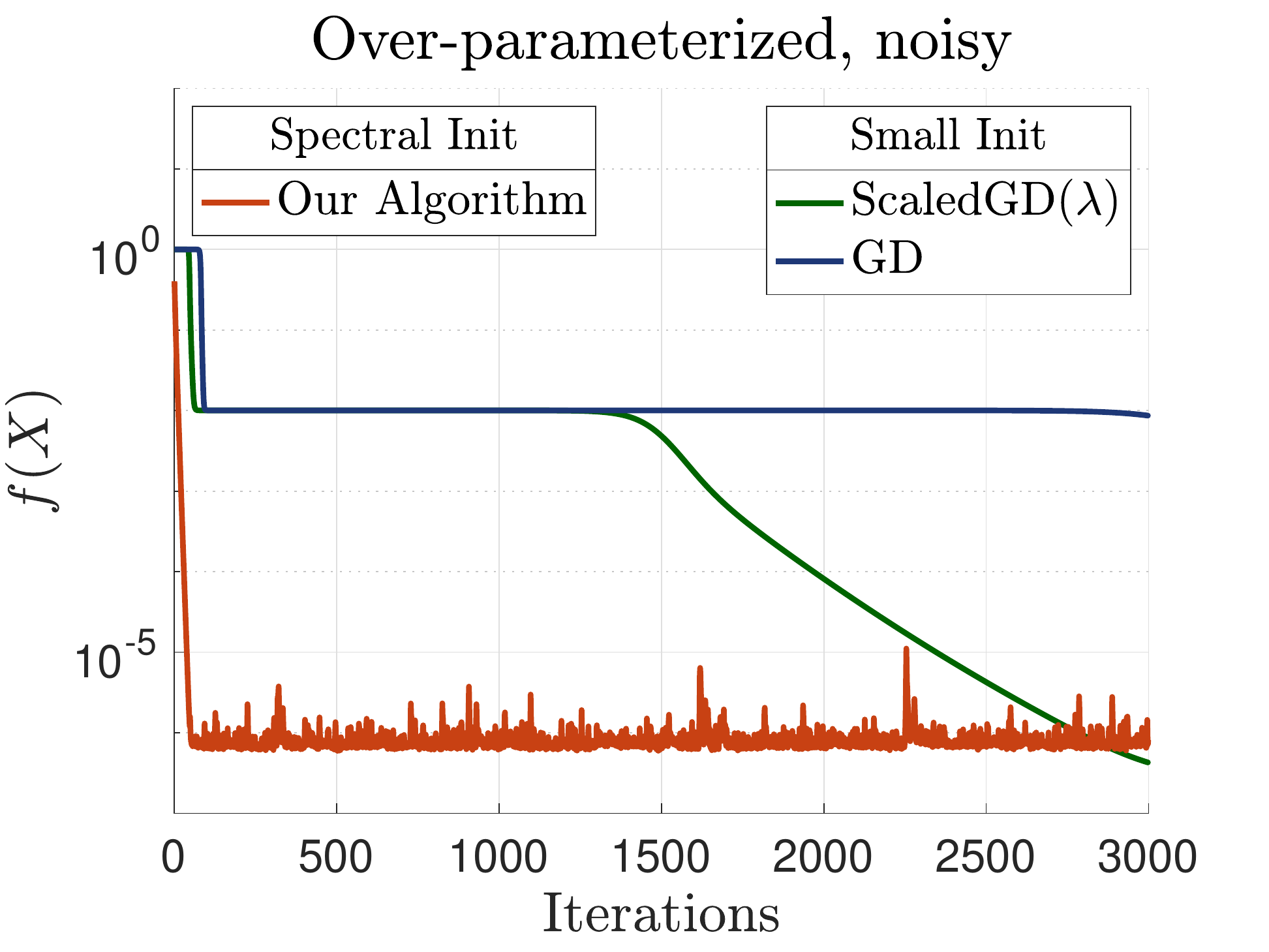}
    \vspace{-2em}
	\caption{\textbf{Convergence of our algorithm (spectral init.), ScaledGD$(\lambda)$ (small init.) and GD (small init.) for Gaussian matrix sensing}. Left: Noiseless measurements. Right: Noisy measurements with noise variance $\sigma=10^{-6}$.}
    \vspace{-0.5em}
    \label{fig:gaussian_matrix_sensing3}
\end{figure*}

  To begin, note that  the objective function in \eqref{eq:obj} can be written as 
\begin{align}
f(X) =f_{c}(X)+\frac{\|\varepsilon\|^{2}}{m}-\frac{2}{m}\langle\mathcal{A}(XX^{T}-M^{\star}),\varepsilon\rangle, \label{eq:fc}
\end{align}
where $f_{c}(X)=\frac{1}{m}\|\mathcal{A}(XX^{T}-M^{\star})\|^{2}$ is defined to be the
objective function with clean measurements that are not corrupted
by noise. 

The first component in our proof consists of showing that the iterates of \eqref{alg:new} can be viewed as gradient descent under a change of norm. In particular, let $P = X^TX + \eta_t I$ be a real symmetric,
positive definite $r\times r$ matrix. We define a corresponding
$P$-norm and its dual $P$-norm on $\mathbb{R}^{n\times r}$
as follows 
\begin{align}
\label{eq:dualnorm}
\|X\|_{P}  \eqdef\|XP^{1/2}\|, \quad\|X\|_{P*} & \eqdef\|XP^{-1/2}\|.
\end{align}
 Consider a descent direction $D$. 
Suppose that the following inequality holds with some constant $L$:
\begin{align}
\label{eq:lip}
f_c(X-\alpha D)  \leq f_c(X)-\alpha\langle\nabla f_c(X),D\rangle+\frac{\alpha^{2}L}{2}\|D\|_{P}^2
\end{align}
and that the PL-inequality holds under the $P$-norm:
$
\|\nabla f_c(X)\|_{P*}^{2}\geq\mu (f_c(X)) \label{eq:gd-P}
$
with $\mu>0$. Then plugging in the descent direction $D = \nabla f_c(X)P^{-1}$
yields linear convergence since 
$
f_c(X-\alpha D) \leq  \left(1- \frac{\mu}{2L}\right) f_c(X). 
$

Therefore, to complete our proof, we need to demonstrate the following conditions:
\begin{enumerate}
	\item The inequality \eqref{eq:lip} holds with some constant $L$
	\item The PL-inequality holds under the $P$-norm: $\|\nabla f_c(X)\|_{P*}^{2}\geq\mu_(f_c(X))$.
\end{enumerate}


First, we consider an ideal case: suppose that there exists some constant $C>1$ such that $\eta_t \leq \sqrt{f_c(X_t)} \leq C\eta_t$. This is exactly the regime for $\eta_t$ where our algorithm is well-behaved: both Lipschitz gradients and the PL-inequality is satisfied in the $P$-norm. The proof of these two facts, especially the second one, is quite involved, but it is similar to the proof of Corollary 5 in \cite{zhang2021preconditioned}, so they are deferred to the appendix.


As a result, in the ideal case where $\eta_t \leq \sqrt{f_c(X_t)} \leq C\eta_t$, if we go in the direction $D = \nabla f_c(X)P^{-1}$, then linear convergence is already achieved.  However, due to noise, the descent direction is $\nabla f(X)$, instead of $\nabla f_c(X)$, since we cannot access the true gradient. Fortunately, if the norm of the gradient is large compared to a statistical error, we can prove that the difference between $\nabla f_c(X)$ and $\nabla f(X)$ is negligible, and our algorithm will still make enough progress at each iteration to ensure linear convergence. 



Essentially, if $\eta_t \approx \sqrt{f_c(X_t)}$, then our algorithm will converge linearly up to some statistical error. 
Unfortunately, $\eta_t \approx \sqrt{f_c(X_t)}$ does not always hold, because \textit{both} $\eta_t$ and $\sqrt{f_c(X_t)}$ are changing. Instead, we have to consider scenarios where $\eta_t$ deviates from this ideal range. To complete the proof of Theorem \ref{thm:main}, we need to show that $\eta_t$ never deviates too far from the ideal range. This is the key difficulty in our proof. Intuitively, if $\sqrt{f_c(X_{t+1})}$ becomes too small compared to $\eta_t$, our algorithm will start to behave more like gradient descent and slow down. Hence with a fixed decay rate, $\eta_t$ will quickly be on the same order as $\sqrt{f_c(X_{t+1})}$ again. As a result, linear convergence is always maintained.

\section{NUMERICAL SIMULATIONS}\label{sec:experiment}
In this section, we compare our method against two state of the art preconditioned methods: PrecGD \cite{zhang2021preconditioned} and ScaledGD($\lambda$) \cite{xu2023power}. To validate our theoretical results, we first perform experiments using Gaussian measurements on a synthetic low-rank matrix. In addition, we also perform experiments on a real-world medical imaging application, specifically in denoising an ultrafast ultrasound scan as shown in Figure~\ref{fig:ultrasound}. We leave the details of this experiment to the appendix.

With spectral initialization, we show that our method is the only algorithm that is able to consistently achieve minimax error at a linear rate. In Figures \ref{fig:gaussian_matrix_sensing1} and \ref{fig:gaussian_matrix_sensing2}, both PrecGD and ScaledGD$(\lambda)$  have trouble converging to the minimax optimal unless a perfect regularization parameter is chosen. Without the perfect choice, they either stagnate at a high noise level, or even diverge. 

In Figure \ref{fig:gaussian_matrix_sensing3}, we see that with small initialization, both ScaledGD$(\lambda)$ and GD seem to stagnate for a significantly long time before converging to the next eigenvalue. Therefore, in cases where a good heuristic initialization is available, throwing such a initialization away and using small initialization instead can come at a great cost, since it can take many iterations to get it back.

\subsection{Gaussian matrix sensing}
In this experiment, we consider a matrix recovery problem on a $10\times 10$ ground truth matrix $M^\star$ with truth rank $r^\star=2$. The condition number of $M^\star$ is set to $\kappa=10^2$. We take measurements on $M^\star$ using linearly independent measurement matrices $A_1,\ldots,A_m$ drawn from the standard Gaussian distribution. In Figure~\ref{fig:gaussian_matrix_sensing1}, \ref{fig:gaussian_matrix_sensing2} and \ref{fig:gaussian_matrix_sensing3}, we set the search rank to be $r=8$ and draw $m=2nr$ measurements from $M^\star$. In noisy setting, we corrupt the measurements with noise $\varepsilon \sim \mathcal{N}(0,\sigma^2)$ where $\sigma=10^{-6}$. For our algorithm, we set $\eta_0=\sqrt{f(X_0)}$ and the decay rate for $\eta_t$ as $\beta = 0.85$ in noiseless case, and $\beta = 0.5$ in noisy case. 
\paragraph{Our algorithm v.s. ScaledGD$(\lambda)$}
In Figure~\ref{fig:gaussian_matrix_sensing1}, we plot the convergence of our algorithm and ScaledGD$(\lambda)$ under four values of $\lambda=\{10^{-2},10^{-4},10^{-6},10^{-9}\}$. We spectrally initialized both methods at the same initial point and set the learning rate to be $\alpha=0.1$. In noiseless case, we see that our algorithm converges to minimax error at a linear rate, while ScaledGD$(\lambda)$ converges to an error of $O(\lambda)$. It is important to note that we cannot set $\lambda$ to be too small because it would cause ScaledGD$(\lambda)$ to diverge or become numerically unstable, as depicted in Figure~\ref{fig:gaussian_matrix_sensing1}. In noisy case, we obtain similar results as in the noiseless case. The only difference is that ScaledGD$(\lambda)$ can converge to the same error as our algorithm when $\lambda=\{10^{-6},10^{-9}\}$ as the minimax error in the noisy case is around $10^{-6}$. However, we again observe that ScaledGD$(\lambda)$ becomes numerically unstable when $\lambda$ is too small.

\paragraph{Our algorithm v.s. PrecGD}
In Figure~\ref{fig:gaussian_matrix_sensing2}, we plot the convergence of our algorithm and PrecGD. We spectrally initialize both methods at the same initial point and set the learning rate to be $\alpha=0.1$. Here, PrecGD is implemented with a proxy variance $\hat{\sigma}$ so that $\eta_t = \sqrt{|f(X_t)-\hat\sigma^2|}$. We see that our algorithm converges linearly to the minimax error in both the noiseless and noisy case. While PrecGD also converges linearly to the minimax error in the noiseless case, in the noisy case, however, its error depends crucially on the value of proxy variance $\hat{\sigma}$; it requires $\hat{\sigma}\approx \sigma$ to achieve minimax error.

\paragraph{Small init. v.s. spectral init.}
In Figure~\ref{fig:gaussian_matrix_sensing3}, we plot the convergence of our algorithm, ScaledGD$(\lambda)$ and GD. In this experiment, our algorithm is initialized using spectral initialization, and both ScaledGD$(\lambda)$ and GD are initialized using small initialization with initialization scale $10^{-12}$. We set the learning rate for all three methods to be $\alpha=0.1$. Again, our algorithm converges linearly to the minimax error in both noiseless and noisy setting. Both ScaledGD$(\lambda)$ and GD learn the solution incrementally (see \cite{jin2023understanding} for a precise characterization of incremental learning) and hence reach the minimax error a lot slower than our algorithm.


\section*{ACKNOWLEDGMENTS}
The authors are grateful to Pengfei Song and YiRang Shin for their help and advice on the ultrafast ultrasound application, and for sharing the rat brain datasets used in our medical imaging experiments. 

Financial support for this work was provided in part by NSF CAREER Award ECCS-2047462.

\bibliography{reference}

\appendix
\onecolumn
\section{Proof of Main Results}
\subsection{Prelminaries}
In addition to the notation used in the main paper, we define some additional notation that will be used throughout the appendix. Let $X$ by an $n\times r$ matrix, and let $P\succ 0$ be a fixed $r\times r$ positive definite matrix. 
We define a corresponding
$P$-norm and its dual $P$-norm on $\mathbb{R}^{n\times r}$
as follows 
\begin{align}
\label{eq:dualnorm}
\|X\|_{P}  \eqdef\|XP^{1/2}\|, \quad\|X\|_{P*} & \eqdef\|XP^{-1/2}\|.
\end{align}
We use $\mathrm{vec}(X)$ to denote the vectorization operator that stacks the column of $X$ into a single column vector. As before, we use $\otimes$ to denote the Kronecker product between two matrices. 
For a scalar-valued function of a matrix, $f(X)$, we use $\nabla^2f(X)[V]$ to denote the Hessian vector product, defined by 
\[
\nabla^2f(X)[V] = \lim_{t\to 0} \frac{\nabla f(X+tV)- \nabla f(X)}{t}. 
\]
Note that here $\nabla^2f(X)[V]$ is a matrix of the same size as $\nabla f(X)$. With a slight abuse of notation, we use lower case letters $x_t$ to denote the vectorized version of $X_t$, so $x_t = \mathrm{vec}(X_t)$. We denote the corresponding gradient by $\nabla f(x_t)$. 

 For symmetric matrix sensing, we denote our ground truth by $M^{\star}=ZZ^{T}\in\R^{n\times n}$. We denote the true rank by $r^* = \mathrm{rank}(M^*)$. Our goal is to recovery $M^\star$ from a small number of measurements of the form $y=\mathcal{A}(ZZ^{T})+\epsilon\in\R^{m}$. Here $\epsilon$ is a vector with independent Gaussian entries with zero mean and variance $\sigma^{2}$. To do so, we minimize the non-convex objective function
\[
f(X)=\frac{1}{m}\|\mathcal{A}(XX^{T})-y\|^{2}=f_{c}(X)+\frac{1}{m}\|\epsilon\|^{2}-\frac{2}{m}\langle\mathcal{A}(XX^{T}-M^{\star}),\epsilon\rangle,
\]
where 
$f_{c}(X) \overset{def}{=} \frac{1}{m}\|\mathcal{A}(XX^{T}-M^{\star})\|^{2}$
is the
objective function with clean measurements that are not corrupted
with noise. Here $X$ is a matrix of size $n\times r$, where $r$ is known as the search rank. 


We make a few additional simplifications on notations. As before, we will use
$\alpha$ to denote the step-size and $D$ to denote the local search
direction. In our proof below, it will often be easier to use the \textit{vectorized} form of our gradient updates: vectorizing both sides of the update $X_{t+1} = X_t - \alpha \nabla f(X_t) (X_t^TX_t + \eta_tI)^{-1}$, we get 
\begin{align*}
\vv{X_{t+1}} & = \vv{X_t}-\alpha \vv{\nabla f(X_t) (X_t^TX_t + \eta_tI)^{-1}} \\
& = \vv X_t - ((X_t^TX_t + \eta_tI) \otimes I)^{-1} \vv(\nabla f(X_t)),
\end{align*}
where in the second line we used the standard identity $\vv(AXB) = (B^T \otimes A) \vv(X)$ for the Kronecker product. Using lower case letters $x$ and $d$ to refer to
$\vv{(X)}$ and $\vv{(D)}$ respectively, the update above can be written as $x_{t+1} = x_t - \alpha \P^{-1} \nabla f(x_t)$, with $\P =  (X_t^TX_t + \eta_tI) \otimes I$.

\subsection{Auxiliary Results}
In this section we collect two results from \cite{zhang2021preconditioned} that will be used in the proof of our main result. The first theorem shows that when the regularization $\eta$ is small, the PL-inequality holds within a small neighborhood around the ground truth. 

\begin{theorem}[Noiseless gradient dominance]
\label{thm:apl}Let $\min_{X}f(X)=0$ for $M^{\star}\ne0$. Suppose
that $X$ satisfies $f(X)\le\rho^{2}\cdot(1-\delta)\lambda_{r^{\star}}^{2}(M^{\star})$
with radius $\rho>0$ that satisfies $\rho^{2}/(1-\rho^{2})\le(1-\delta^{2})/2$.
Then, we have
\[
\eta\le C_{ub}\|XX^{T}- M^\star \|_{F}\quad\implies\quad\|\nabla f(X)\|_{P*}^{2}\ge2\mu f(X)
\]
where
\begin{equation}
\mu =\left(\sqrt{\frac{1+\delta^{2}}{2}}-\delta\right)^{2}\cdot\min\left\{ \left(\frac{C_{1}}{\sqrt{2}-1}\right)^{-1},\left(1+3C_{1}\sqrt{\frac{(r-r^{\star})}{1-\delta^{2}}}\right)^{-1}\right\}. \label{eq:mueq}
\end{equation}
Here $C_{1}$ is a constant that only depend on $\delta$. 
\end{theorem}
We recall that in the theorem above $\|\cdot\|_{P^*}$ denotes the dual norm of $\|\cdot\|_P$ defined in \eqref{eq:dualnorm}, with $P = X^TX+\eta I$. Essentially, this theorem says that when $\eta$ is small compared to the true error $\|XX^T-M^\star\|$, the PL-inequality is restored under the local norm defined by $P$. The difficulty of applying this theorem directly in our case arises from two issues: first, if $\eta$ is too small, then the gradients of $f(X)$ are no longer Lipschitz under the $P$-norm. As a result, the iterates can diverge. Moreover, it is very difficult to gauge the `right' size of $\eta$ in the noisy setting, since we have no access to the true error. The proof of Theorem \ref{thm:apl} can be found in \cite{zhang2021preconditioned} so we do not repeat it here.

We also state an lemma from \cite{zhang2021preconditioned} that directly characterizes the progress of gradient descent at each iteration in a fashion similar to the descent lemma (see e.g. \cite{nesterov2018lectures}). For general smooth functions with Lipschitz gradients, the decrement in the function value at each iteration can be characterized by a quadratic upper bound (the so-called descent lemma). However, for matrix sensing, we can in fact obtain a tighter upper bound because $f_{c}(X-\alpha D)$ itself is just a quartic polynomial. This allows us to characterize the progress made at each iteration directly, using the following result. 
\begin{lemma}
\label{lem:descent} For any descent direction $D\in\mathbb{R}^{n\times r}$
and step-size $\alpha>0$ we have 
\begin{align}
f_{c}(X-\alpha D)&\le f_{c}(X)-\alpha \inp{\nabla f_{c}(X)}{D}+\frac{\alpha^{2}}{2}\inp{D}{\nabla^{2}f_{c}(X)[D]}\\
&+\frac{(1+\delta)\alpha^{3}}{m}\|D\|_{F}^{2}\left(2\|DX^{T}+XD^{T}\|_{F}+\alpha\|D\|_{F}^{2}\right).\label{eq:f_descent}
\end{align}
\end{lemma}

The proof of this lemma is quite straightforward so we do not repeat it here. The Lemma follows simply from expanding the function $f_c(X-\alpha D)$ and bounding the third and fourth order terms using the restricted isometry property of $\mathcal{A}$.

In section 3 of the main paper, we sketched out the main idea behind our proof of Theorem \ref{thm:main}. In particular, we stated that our proof mainly consists of two parts. First, in the ideal case where $\eta_t \leq \sqrt{f_c(X_t)} \leq C\eta_t$, if we go in the direction $D = \nabla f_c(X)P^{-1}$, then linear convergence is already achieved in the sense that the function value decreases by a constant factor. Second, in the case where $\eta_t \leq \sqrt{f_c(X_t)} \leq C\eta_t$ is no longer satisfied, we want to show that $\eta_t$ will automatically return to the ideal interval $\eta_t \leq \sqrt{f_c(X_t)} \leq C\eta_t$ after a few iterations. Thus overall, we have linear convergence. In the following sections, we make these two parts of our proof precise.

\subsection{Linear Convergence in Ideal Case}
For the PrecGD algorithm of \cite{zhang2021preconditioned} to succeed, it is crucial that the regularization parameter satisfies  $\eta_t \leq \sqrt{f_c(X_t)} \leq C\eta_t$. If so, then within a local neighborhood of the ground truth, Theorem \ref{thm:apl} and Theorem \ref{lem:descent} can be used to establish linear convergence in the noiseless setting. However, as we have argued in the main paper, this requirement for $\eta_t$ is difficult if not impossible to maintain explicitly in the noisy setting. This makes it extremely difficult for PrecGD to achieve a minimax optimal error.

Our main result, Theorem \ref{thm:main}, states that  letting $\eta_t$ decay with some constant rate $\beta$ suffices to guarantee the linear convergence of our algorithm, even in the noisy setting. One of our key observations is that the condition  $\eta_t \leq \sqrt{f_c(X_t)} \leq C\eta_t$ \textit{does not} have to be satisfied at all times. In fact, we can allow $\eta_t$ to dip below $\sqrt{f_c(X_t)}$ in our algorithm, because of the ``coupling'' effect that we discussed previously: $\eta_t$ can never deviate too far from $\sqrt{f_c(X_t)}$. 

Therefore, in our proof of Theorem \ref{thm:main}, we will consider two cases: 
\begin{enumerate}
	\item $\eta_t \leq \sqrt{f_c(X_t)} \leq C\eta_t$ 
	\item $\sqrt{f_c(X_t)} \leq \eta_t$.
\end{enumerate}
The first case is the ``good'' situation, because the conditions for linear convergence is satisfied by assumption. In this case, we show that as long as the gradient is large compared to the noise, i.e., $
\|\nabla f_c(X_t)\|_P^* \gtrsim \sqrt{\frac{\sigma^{2}rn\log n}{m}},
$  
our algorithm will converge linearly. This behavior is stated rigorously in the following lemma.

\begin{lemma}
\label{lem:noise}
Suppose that at the $t$-th iteration, the regularization parameter $\eta_t$ satisfies $\eta_t \leq \sqrt{f_c(X_t)} \leq C\eta_t$ for some $C>1$. Furthermore, suppose that 
$
\|\nabla f_c(X_t)\|_P^* \gtrsim \sqrt{\frac{\sigma^{2}rn\log n}{m}}.
$
Then for $\alpha \leq 1/L$, with high probability we have 
\begin{align*}
f_{c}(X_{t+1}) 
& \leq \left(1-\frac{\mu}{2L}\right)f_{c}(X_t).
\end{align*}
Here $L = O(C^2)$ is the constant defined in \eqref{eq:lp}, and $\mu$ is the constant defined in Theorem \ref{thm:apl}.
\end{lemma}

The proof of Lemma \ref{lem:noise} is long but it is mainly computational. The overall idea is similar to the proof of Theorem 20 in \cite{zhang2021preconditioned}: our goal is to show that when the local norm of the gradient is large compared to the noise level, the decrement we make at each iteration `overcomes' the error caused by the noisy measurements. Our main tool here is Lemma \ref{lem:descent}, which allows us to directly compute the decrement and bound the error terms. 


It turns out that in this proof, it will be slightly easier to deal with the vectorized version of this problem: we use $f(x)$ to denote original objective function $f(X)$ as a function of the vector $x = \vv(X)$. 
Consequently, we write $f(x)\in\R^{nr}$ and $\nabla f(x)\in R^{nr}$
as the vectorized versions of $f(X)$ and its gradient. We use the same vectorized notation for the ``true'' function value $f_{c}(X)$. Thus, in vectorized form, the iterates of our algorithm can be written as
\[
x_{k+1} = x_k - \alpha  \P^{-1} \nabla f(x), \quad \text{where } \P=(X^{T}X+\eta I_{r})\otimes I_{n}.
\]
We note that all the norms we consider remain unchanged after vectorization, meaning that  
$
\|\nabla f(x)\|_{P}=\|\nabla f(X)\|_{P}$ and $\|\nabla f(x)\|_{P^{*}}=\|\nabla f(X)\|_{P^{*}}.
$
Now we are ready to prove this lemma.

\begin{proof}
The main idea of the proof is to use the inequality $\eta_t \leq \sqrt{f_c(X_t)} \leq C\eta_t$ to bound the progress of our algorithm at each iteration. In particular,  when $\eta_t$ is small, i.e., $\eta_t \leq \sqrt{f_c(X_t)}$, then Theorem \ref{thm:apl} guarantees the gradient dominance. On the other hand, the lower bound $\sqrt{f_c(X_t)} \leq C\eta_t$ allows us to apply Lemma \ref{lem:descent} to guarantee that the step-size $\alpha$ can be large enough so that we get linear convergence. 

First, note that \textit{vectorized} version of the gradient update $X_{+} = X -\alpha D$ (where $D = \nabla f(X)P^{-1}$) 
can be
written as $x_{+}=x-\alpha d$, where 
\begin{equation}
\begin{aligned}d & =\vv{(\nabla f(X)P^{-1})} &=\P^{-1}\nabla f_{c}(x)-\frac{2}{m}\P^{-1}\left(I_{r}\otimes\sum_{i=1}^{m}\epsilon_{i}A_{i}\right)x.
\end{aligned}
\label{eq:d_gd}
\end{equation}
Here we have dissected the gradient descent direction into two parts: $\P^{-1}\nabla f_{c}(x)$, which corresponds to ``correct'' gradient and a remaining error term $\P^{-1}\mathcal{E}(x) $, where
\[
\mathcal{E}(x) \overset{def}{=} \frac{2}{m}\left(I_{r}\otimes\sum_{i=1}^{m}\epsilon_{i}A_{i}\right)x.
\]
In other words we have $d = \P^{-1}(\nabla f_c(x)-\mathcal{E}(x))$. If $\mathcal{E}(x) = 0$, then our proof reduces to the noiseless case. Here we want to show that the error is small compared the decrement we make in the function value. As we will see, this happens precisely in the regime where the gradient is large, i.e.,  $\|\nabla f_c(X)\|_P^* \gtrsim \sqrt{\frac{\sigma^{2}rn\log n}{m}}.$

In vectorized notation, Lemma \ref{lem:descent} can be written as 
\begin{equation}
f_{c}(x-\alpha d)\le f_{c}(x)-\alpha\nabla f_{c}(x)^{T}d+\frac{\alpha^{2}}{2}d^{T}\nabla^{2}f_{c}(x)d+\frac{(1+\delta)\alpha^{3}}{m}\|d\|^{2}\left(2\|\J d\|+\alpha\|d\|^{2}\right),\label{eq:noisy_descent_vec}
\end{equation}
where we define $\J:\R^{nr}\to\R^{n^{2}}$ as the linear operator
satisfying $\J d=\vv(XD^{T}+DX^{T})$ (recall that $d = \vv(D)$). Now setting $d = \P^{-1}(\nabla f_c(x)-\mathcal{E}(x))$ in the formula above yields 
\begin{align*}
f_{c}(x-\alpha d)\leq & f_{c}(x)-\alpha\|\nabla f_{c}(x)\|_{P*}^{2}+T_{1}+T_{2}+T_{3}
\end{align*}
where 
\begin{align*}
T_{1}= & \alpha\nabla f_{c}(x)^{T}\P^{-1}\E(x)\\
T_{2}= & \frac{\alpha^{2}}{2}\Big(\nabla f_{c}(x)^{T}\P^{-1}\nabla^{2}f_{c}(x)\P^{-1}\nabla f_{c}(x)+\E(x)^{T}\P^{-1}\nabla^{2}f_{c}(x)\P^{-1}\E(x)\\
 & -2\nabla f_{c}(x)^{T}\P^{-1}\nabla^{2}f_{c}(x)\P^{-1}\E(x)\Big)\\
T_{3}= & (1+\delta)\alpha^{3}\left(\|\P^{-1}\nabla f_{c}(x)-\P^{-1}\E(x)\|^{2}\right)\left(2\|\J\P^{-1}\nabla f_{c}(x)\|+2\|\J\P^{-1}\E(x)\|\right.\\
 & +\left.\alpha\|\P^{-1}\nabla f_{c}(x)-\P^{-1}\E(x)\|^{2}\right).
\end{align*}
Our goal is to show that all three terms $T_1,T_2,T_3$ are small compared to the decrement that we make at each iteration. 
 The key observation here is that all of these terms depend on $\E(x)$ and $\P$. With the right choice of $\eta$, i.e., with $\sqrt{f_c(X_t)} \leq C\eta_t$, the preconditioner $\P$ is well-conditioned, so that all the errors in $T_1, T_2, T_3$ will remain small as long as $\E(x)$ is small. Specifically, we can bound the error term as 
\begin{align*}
\|\E(x)\|_{P^{*}}^{2} & =\E(x)^{T}\P^{-1}\E(x)=\left\Vert \left(\frac{2}{m}\sum_{i=1}^{m}\epsilon_{i}A_{i}\right)X(X^{T}X+\eta I)^{-1/2}\right\Vert _{F}^{2}\\
 & \leq\left\Vert \left(\frac{2}{m}\sum_{i=1}^{m}\epsilon_{i}A_{i}\right)\right\Vert _{2}^{2}\left\Vert X(X^{T}X+\eta I)^{-1/2}\right\Vert _{F}^{2}\\
 & \overset{(i)}{}\leq C_e\frac{\sigma^{2}n\log n}{m}\left(\sum_{i=1}^{r}\frac{\sigma_{i}^{2}(X)}{\sigma_{i}(X)^{2}+\eta}\right)\\
 & \leq C_e\frac{\sigma^{2}rn\log n}{m},
\end{align*}
where $C_e$ is an absolute constant and (i) follows from a standard concentration bound (see \cite{candes2011tight} or Lemma 16 of \cite{zhang2021preconditioned}). 

Now,  denoting $\Delta=\|\nabla f_{c}(x)\|_{P*}$ and using the bound for the error above, we get after some computations that 
\begin{align*}
 & T_{1}\leq\alpha\Delta\sqrt{\frac{C_e\sigma^{2}rn\log n}{m},}\\
 & T_{2}\leq2\alpha^{2}L_{\delta}\Delta^{2}+2\alpha^{2}L_{\delta}\frac{\sigma^{2}rn\log n}{m}\\
 & T_{3}\leq\frac{4(1+\delta)\alpha^{3}}{\eta}\left(\Delta^{2}+\frac{C_e\sigma^{2}rn\log n}{m}\right)\left(\frac{\alpha\Delta^{2}}{\eta}+\frac{\alpha C_e\sigma^{2}rn\log n}{\eta m}+2\sqrt{2}\Delta+2\sqrt{2}\sqrt{\frac{C_e\sigma^{2}rn\log n}{m}}\right).
\end{align*}
Here $L_\delta$ is a constant that depends only on the RIP constant $\delta$. Now plugging these error bounds back into (11) yields
\begin{align}
f_{c}(x-\alpha d) & \leq f_{c}(x)-\alpha\Delta^{2}+\alpha\Delta\sqrt{\frac{C\sigma^{2}rn\log n}{m}}+2\alpha^{2}L_{\delta}\Delta^{2}+2C\alpha^{2}L_{\delta}\frac{\sigma^{2}rn\log n}{m}\nonumber \\
 & +\frac{4(1+\delta)\alpha^{3}}{\eta}\left(\Delta^{2}+\frac{C\sigma^{2}rn\log n}{m}\right)\left(\frac{\alpha\Delta^{2}}{\eta}+\frac{\alpha C\sigma^{2}rn\log n}{\eta m}+2\sqrt{2}\Delta+2\sqrt{2}\sqrt{\frac{C\sigma^{2}rn\log n}{m}}\right).\label{eq:eq_f2}
\end{align}

In the case $\Delta\geq2\sqrt{\frac{C_e \sigma^{2}rn\log n}{m}}$ all the terms above can be bounded so that the decrement in the function value dominates all the error. In particular, plugging this lower bound into the inequality above yields
\begin{align*}
f_{c}(x-\alpha d)\leq f_{c}(x)-\frac{\alpha}{2}\Delta^{2}\left(1-\frac{5}{2}L_{\delta}\alpha-60\sqrt{2}\alpha^{2}(1+\delta)-25\alpha^{3}(1+\delta)^{2}\right).
\end{align*}
Now, assuming that the step-size satisfies 
\begin{equation}
\label{eq:lp}
\alpha\leq\min\left\{ \frac{L_{\delta}}{60\sqrt{2}(1+\delta)+25(1+\delta)^{2}},\frac{1}{7L_{\delta}}\right\}  \overset{def}{=} \frac{1}{L}
\end{equation}
we obtain
$
f_{c}(x-\alpha d)\leq f_{c}(x)-\frac{t\Delta^{2}}{4}\leq\left(1-\frac{\alpha\mu}{4}\right)f_{c}(x),
$
where in the last step we used the fact that $\eta_t \leq \sqrt{f_c(X_t)}$, so the conditions of Theorem \ref{thm:apl}  are satisfied, so gradient dominance holds. This completes the proof. 
	
\end{proof}

\section{Proof of Theorem \ref{thm:main}}
In this section we provide a complete proof of Theorem \ref{thm:main}, filling out some of the missing details left out in the main paper. 
\begin{proof}
Let $T>0$ be the smallest index such that $\eta_T < 2 \sqrt{\frac{C_0 \sigma^{2}rn\log n}{\mu m}}$. Suppose that $t<T$.
Similar to the noiseless case, we will show that there exists some constant $C>1$, which depends only on $\delta$, such that the following holds: if at the $t$-th iterate we have $\sqrt{f_c(X_t)} \leq C\eta_t$, then $\sqrt{f_c(X_{t+1})} \leq C\eta_{t+1}$. As before, this implies that $f_c(X_t) \leq C^2\beta^{2t}\eta_0$ for all $t\leq T$. 

At the $t$-iterate, suppose that $\sqrt{f_c(X_t)} \leq C\eta_t$. We consider two cases: 
\begin{enumerate}
	\item $\eta_t \leq \sqrt{f_c(X_t)} \leq C\eta_t$ 
	\item $\sqrt{f_c(X_t)} \leq \eta_t$.
\end{enumerate}
 We will show that in either case, the next iterate satisfies $\sqrt{f_c(X_{t+1})} \leq C\eta_{t+1}.$ 
 The core idea behind this proof is that the values of $\eta_t$ and $\sqrt{f_c(X_{t+1})}$ are ``coupled'', meaning that they can not deviate too far from each other. In the first case, where $\eta_t \leq \sqrt{f_c(X_t)} \leq C\eta_t$, the behavior of our algorithm plus is exactly the same as PrecGD, since $\eta_t$ is bounded both above and below by a constant factor of $\sqrt{f_c(X_{t+1})}$. Thus, according to Lemma \ref{lem:noise}, we converge linearly (at least for the current iteration). In fact, we have chosen the decay rate $\beta$ so that $f_c(X)$ will decay faster than $\beta$ when $\eta_t \leq \sqrt{f_c(X_t)} \leq C\eta_t$. Specifically, we have 
\[
\sqrt{f_{c}(X_{t+1})} \leq \sqrt{\left(1-\frac{\mu}{4L}\right)} \sqrt{f_{c}(X_t)} \leq \beta \cdot C\eta_t  = C\eta_{t+1}
\]
where the second inequality follows from $\beta = \sqrt{\left(1-\frac{\mu}{8L}\right)}$ and the assumption $\sqrt{f_{c}(X_t)} \leq C\eta_t$. Thus, in this case, 
$\sqrt{f_{c}(X_{t+1})}$ will continue to be upper bounded by $C\eta_{t+1}$. If this remains true for all $t$, then we are already done since $\eta_{t}$ decays exponentially, which means that the function value will also decay exponentially fast. However, if  $\sqrt{f_{c}(X_{t+1})}$ decays too fast, the condition for applying Lemma \ref{lem:noise}, i.e., $\eta_t \leq \sqrt{f_c(X_t)} \leq C\eta_t$, will no longer hold. However, in this case, the function values are still decaying monotonically. Since the stepsize satisifies $\alpha \leq 1/L$, where $L$ is the constant defined in \eqref{eq:lp}, we can use Lemma \ref{lem:noise} again to get $f_c(X_{t+1}) \leq f_c(X_t)$. Thus   
\[
\sqrt{f_c(X_{t+1})} \leq \sqrt{f_c(X_{t})} \leq \eta_t = \beta^{-1} \eta_{t+1} \leq C \eta_{t+1}. 
\]
Here we note that for the last step to hold we need 
$
\beta^{-1}<C$, which is equivalent to $\sqrt{1- \frac{\mu}{4L}} \cdot C > 1$. In fact, this is the key step that keeps us from choosing the decay rate $\beta$ to be too small so that we get any linear convergence rate we like. By definition $L  = C_\delta \cdot C^2$, where $C_\delta$ is a constant that only depends on $\delta$. Thus this condition is always satisfied for some $C\geq C_{lb}$, where $C_{lb}$ is a constant lower bound that only depends on $\delta$. 

Finally, at the $T$-th iteration, we have $\sqrt{f_c(X_T)} \leq C\cdot \eta_T \lesssim  \sqrt{\frac{\sigma^{2}rn\log n}{m}}$. Now for all $t>T$, we again consider two cases: 
\begin{enumerate}
	\item $\|\nabla f_c(X_t)\|_P^* \leq \sqrt{\frac{\sigma^{2}rn\log n}{m}}$ 
	\item $\|\nabla f_c(X_t)\|_P^* \geq \sqrt{\frac{\sigma^{2}rn\log n}{m}}$.
\end{enumerate}
In the first case, we can use Theorem \ref{thm:apl} to conclude that $\mu f(X_t) \leq \left(\|\nabla f_c(X_t)\|_P^*\right)^2$. Since $\mu$ is a constant, we have $f(X_t) \lesssim \frac{\sigma^{2}rn\log n}{m}$. Now consider second case. Here we can apply Lemma \ref{lem:noise} again which guarantees that $f(X_{t+1}) \leq f(X_t)$, so the function value is decreasing. Consequently, we have $f(X_t) \lesssim \frac{\sigma^{2}rn\log n}{m}$ for all $t>T$. This completes the proof.


\end{proof}
\section{Experimental details}\label{sec:dataset}
\paragraph{Experimental setups} We perform all the experiments in this paper on an Apple MacBook Pro, running a silicon M1 pro chip with 10-core CPU, 16-core GPU, and 32GB of RAM. We implement our algorithm in MATLAB R2021a.
\paragraph{Initialization}
\begin{itemize}
    \item \textbf{Spectral initialization:} For spectral initialization with respect to a ground truth matrix $M^\star = Q\Sigma Q^T$, we initialized $X_0=Q\Sigma^{1/2}+0.1\cdot\hat Q$, where $\hat Q\in\mathbb R^{n\times r}$ is drawn from standard Gaussian.
    \item \textbf{Small initialization:} For small initialization, we set $X_0=\hat\alpha\cdot\hat Q$, where $\hat\alpha$ is the initialization scale and $\hat Q$ is a $n\times r$ matrix drawn from standard Gaussian.
\end{itemize}
\paragraph{Datasets}
The datasets we use for the experiments in the main paper and Appendix~\ref{sec:additional_exp} are described below.
\begin{itemize}
    \item \textbf{Gaussian matrix sensing:} For the experiment results shown in Figure~\ref{fig:gaussian_matrix_sensing1}, \ref{fig:gaussian_matrix_sensing2}, \ref{fig:gaussian_matrix_sensing3} and \ref{fig:gaussian_matrix_sensing} we synthetically generate a $10\times 10$ ground truth matrix $M^\star$. The rank of $M^\star$ is set to 2. To generate $M^\star$, we first randomly generate an orthonormal matrix $Q\in\mathbb R^{10\times 2}$ and then set $M^\star = Q\Sigma Q^T$. For Figure~\ref{fig:gaussian_matrix_sensing1}, \ref{fig:gaussian_matrix_sensing2} and \ref{fig:gaussian_matrix_sensing3}, we set $\Sigma=\mathrm{diag}(1,10^{-2})$ so that $M^\star$ is ill-conditioned with condition number $\kappa=10^{2}$. For Figure~\ref{fig:gaussian_matrix_sensing}, we set $\Sigma=\mathrm{diag}(1,1)$ so that $M^\star$ is well-conditioned with condition number $\kappa=1$.
    \item \textbf{1-bit matrix sensing:} For the experiment results shown in  Figure~\ref{fig:1bit_matrix_sensing}, the ground truth matrix $M^\star$ is exactly the same as the one in Figure~\ref{fig:gaussian_matrix_sensing}.
    \item \textbf{Phase retrieval:} For the experiment results shown in Figure~\ref{fig:phase_retrieval}, we synthetically generate a length $10$ complex ground truth vector $z$ and set $M^\star = zz^T$. The real and imaginary parts of $z$ are drawn from standard Gaussian.
    \item \textbf{Ultrafast ultrasound image denoising task:} For the experiment results shown in Figure~\ref{fig:ultrasound} and \ref{fig:ultrasound_sampling}. We take an ultrafast ultrasound scan on a rat brain provided from our collaborator. The ultrasound scan consists of 2400 frames of size $200\times 130$ images. We note that in order to show the entire ultrasound scan in 2D, in Figure~\ref{fig:ultrasound} and \ref{fig:ultrasound_sampling}, the ultrasound scans are shown in power Doppler \cite{bercoff2011ultrafast}. In particular, let $M_i\in\mathbb R^{200\times 130}$ be the $i$-th frame of the ultrasound scan, the power Doppler of the scan is defined as  $P_{M}=20\log\left(\mu_M\sum_{i} M_i^2 \right)$ where $\mu_M$ is a normalization constant that normalizes the entries in $\mu_M\sum_{i} M_i^2$ to between 0 and 1. Here, $M_i^2$ denotes the elementwise squaring. In this case, the ultrasound image denoising problem can be treated as a low-rank matrix completion problem on a size $26000\times 2400$ ground truth matrix $M^\star$, which we will elaborate in Appendix~\ref{app:ultrasound}.
\end{itemize}

\section{Ultrafast ultrasound image denoising task}\label{app:ultrasound}
Ultrafast ultrasound is an advanced imaging technique that leverages high frame rate of plane wave imaging, reaching up to thousands of frames per second. The significant increase in frame rate has revolutionized ultrasound imaging, particularly in ultrafast Doppler \cite{bercoff2011ultrafast}, providing enhanced temporal resolution for precise evaluation of high-speed blood flows and improved sensitivity in detecting subtle flow within small vessels. However, clutter signals originated from stationary and slow moving tissue introduce significant artifacts during the acquisition of ultrasound image, preventing it from capturing a clear visualization of vascular paths, and measuring blood flows in small vessels. Therefore, effective denoising techniques for removing these artifacts are often required to obtain a high-quality ultrasound image, in order to minimize the chances of diagnostic
errors, detect subtle changes or anomalies, and reduce
the need for repeated scans or reanalysis.

One effective denoising technique to suppress clutter signals is based on computing the truncated SVD on a space-time matrix \cite{demene2015spatiotemporal}. In particular, let $M_i$ be the $i$-th frame of the $m$ frames ultrasound scan, \cite{demene2015spatiotemporal} proposed to compute the SVD on the space-time matrix $M=[\mathrm{vec}(M_1)\ldots\mathrm{vec}(M_{m})]=Q\Sigma S^T$ where $Q$ and $S$ are orthonormal matrices and $\Sigma$ is diagonal. The noiseless space-time matrix $M^\star$ can then be computed by keeping the top $r$ singular values in $\Sigma$, i.e. $M^\star=Q_r\Sigma_rS_r^T$ where $Q_r$ and $S_r$ denotes the first $r$ columns of $Q$ and $S$, respectively, and $\Sigma_r$ denotes the top $r\times r$ block of $\Sigma$. This technique is effective as signals from stationary and slow moving tissue generally correspond to low frequency components in the spectral domain, as they change slowly over time. However, due to the high frame rates of the ultrafast ultrasound, with prolonged acquisition, it would require sufficiently large memories to store $M$. In many cases, computing the SVD on a large $M$ also become computationally prohibit as its complexity is cubic in the number of frames $m$. 

One possible way to address these limitations is to downsampled the space-time matrix $M$ and then use it to approximate the noiseless space-time matrix $M^\star$ through low-rank matrix completion. Specifically, we treat the noiseless space-time matrix $M^\star$ as the ground truth, and perform low-rank matrix completion $M^\star= UV^T$ using noisy measurements $y_{ij}=M^\star_{i,j}+\epsilon_{i,j}=M_{i,j}$. Here, $U$ and $V$ are matrices with exactly $r\ll m$ columns. Notice that if the matrix completion problem achieves minimax error, it is exactly coincides with the truncated SVD, i.e. $U=Q_r\Sigma^{1/2}$ and $V=S_r\Sigma^{1/2}$.

In the experimental results shown in Figure~\ref{fig:ultrasound}, we are interested in denoising a 60 megapixel (2400-frames, $200\times 130$ pixels per frame) ultrafast ultrasound image by running 30 iterations of the low-rank denoising procedure in \cite{demene2015spatiotemporal} with 50\% sampling rate. As described above, this ultrasound image denoising task can be viewed as a matrix completion problem on a size $26000\times 2400$ ground truth matrix $M^\star$ given 50\% of its noisy entries. In our experiment, we randomly sample (without replacement) 50\% of the entries in $M=[\mathrm{vec}(M_1)\ldots\mathrm{vec}(M_{2400})]$ as our noisy measurements: each measurement takes the form $y_{ij}=\inp{e_ie_j^T}{M}=\inp{e_ie_j^T}{M^\star}+\epsilon_{i,j}$, which is a noisy measurement on $M^\star$. We approximate the noiseless ground truth by first setting $M^\star= UV^T$ and minimize the following loss function over rank-$r$ matrices $U$ and $V$
\[
f(U,V)=\frac{1}{|\Omega|}\sum_{(i,j)\in\Omega}\left(\inp{e_ie_j^T}{UV^T} - y_{ij}\right)^2
\]
where the set $\Omega=\{(i,j)\}$ contains indices for which we know the value of $M_{ij}$. 

In Figure~\ref{fig:ultrasound}, we show the ultrasound image recovered from our algorithm, PrecGD (best previous) and GD (no preconditioning) in power Doppler \cite{bercoff2011ultrafast}. As shown in Figure~\ref{fig:ultrasound}, our algorithm is the only algorithm that achieves the best denoising effect, making the image even sharper. We also emphasize that the per-iteration cost of our algorithm is almost identical to gradient descent.  All three experiments take approximately 3 minutes.

\paragraph{Time complexity for gradient evaluation} In this experiment, because the measurements is of the form $y_{ij}=\inp{e_ie_j^T}{M}$, the two gradient terms $\nabla_{U} f(U,V)$ and $\nabla_{V} f(U,V)$ in \eqref{alg:1} can be efficiently calculated in $O(n_1r|\Omega|)$ time and $O(n_2r|\Omega|)$ time, respectively. Here, we let $n_1$ denote the number of rows in $U$ and $n_2$ denote the number of rows in $V$. To see why this is the case, observe that the two gradient terms can be expressed as $\nabla_{U} f(U,V) = EV$ and $\nabla_{V} f(U,V) = E^TU$ where
\[
E = \frac{2}{|\Omega|}\sum_{(i,j)\in\Omega} \inp{e_ie_j^T}{UV^T}\cdot e_ie_j^T
\]
is a size $n_1\times n_2$ sparse matrix with exactly $|\Omega|$ nonzero entries, which can be efficiently formed in $O(r|\Omega|)$ time and $O(|\Omega|)$ memory. Hence, despite the large number of measurements in this experiment ($|\Omega|=31.2$ million), in our practical implementation, evaluating both gradient terms at each iteration only takes approximately 6 seconds.

\paragraph{Ultrasound image denoising task}
In the experiment results shown in Figure~\ref{fig:ultrasound}, we set the search rank to be $r=100$, so that $U$ is a size $26000\times 100$ matrix and $V$ is a size $2400\times 100$ matrix. We apply our algorithm, PrecGD and GD to minimize $f(U,V)$ for 30 iterations. Here, our algorithm is implemented with $\beta=0.05$, and PrecGD is implemented with proxy variance $\hat \sigma=5\times 10^{-3}$ so that $\eta_t=\sqrt{|f(U_t,V_t)-\hat \sigma^2|}$. The learning rate for our algorithm, PrecGD and GD are chosen to be as large as possible. For our algorithm and PrecGD, the learning rate is set to be $\alpha=10^7$. For GD, the learning rate is set to be $\alpha=10^3$.

\newpage
\section{Additional experiments}\label{sec:additional_exp}
In this section, we perform an additional experiment on Gaussian matrix sensing. We also perform additional experiments to gauge the performance of our algorithm for applications outside of the assumptions of our theoretical results. In particular, we consider two common problems considered in the existing literature that do not satisfy conditions under which Theorem \ref{thm:main} applies: phase retrieval and 1-bit matrix sensing. For these problems, we see almost identical results to Gaussian matrix sensing: our algorithm succeeds in converging to a minimax optimal error, while GD, PrecGD and ScaledGD$(\lambda)$ struggle. 

\subsection{Gaussian matrix sensing}
The problem formulation is described in the main paper. In this experiment, we take 80 measurements $y_i=\inp{A_i}{M^\star}$ on the ground truth matrix $M^\star\in\mathbb R^{10\times 10}$ using 80 linearly independent measurement matrices $A_i\in\mathbb R^{10\times 10}$ drawn from standard Gaussian. Substituting $M^\star = XX^T$, the loss function for Gaussian matrix sensing is defined as
\[
f(X)=\frac{1}{80}\sum_{i=1}^{80}\left(\inp{A_i}{XX^T}-y_i\right)^2.
\]
We perform Gaussian matrix sensing under four different settings.
\paragraph{The exactly-parameterized, noiseless case} Recall that the truth rank of $M^\star$ is 2. In the exactly-parameterized case, we set $X$ to be a size $10\times 2$ matrix and minimize $f(X)$ using our algorithm, PrecGD, ScaledGD$(\lambda)$ and GD for 500 iterations. We set $\beta=0.1$ in our algorithm, and $\lambda=0$ in ScaledGD$(\lambda)$. The learning rate for all four methods are set to $\alpha=0.1$.

\paragraph{The over-parameterized, noisy case} In this setting, we corrupt the measurements with noise $\varepsilon_i\sim\mathcal N(0,10^{-6})$ such that $y_i=\inp{A_i}{M^\star}+\varepsilon_i$. We set $X$ to be a size $10\times 4$ matrix and minimize $f(X)$ using our algorithm, PrecGD, ScaledGD$(\lambda)$ and GD for 500 iterations. Here, our algorithm is implemented with $\beta=0.1$, PrecGD is implemented with proxy variance $\hat \sigma=10^{-5}$ so that $\eta_t=\sqrt{|f(X_t)-\hat \sigma^2|}$, and ScaledGD$(\lambda)$ is implemented with $\lambda=0.01$. The learning rate for all four methods are set to $\alpha=0.1$.

\begin{figure*}[h]
	\centering
	\includegraphics[width = 0.5\textwidth]{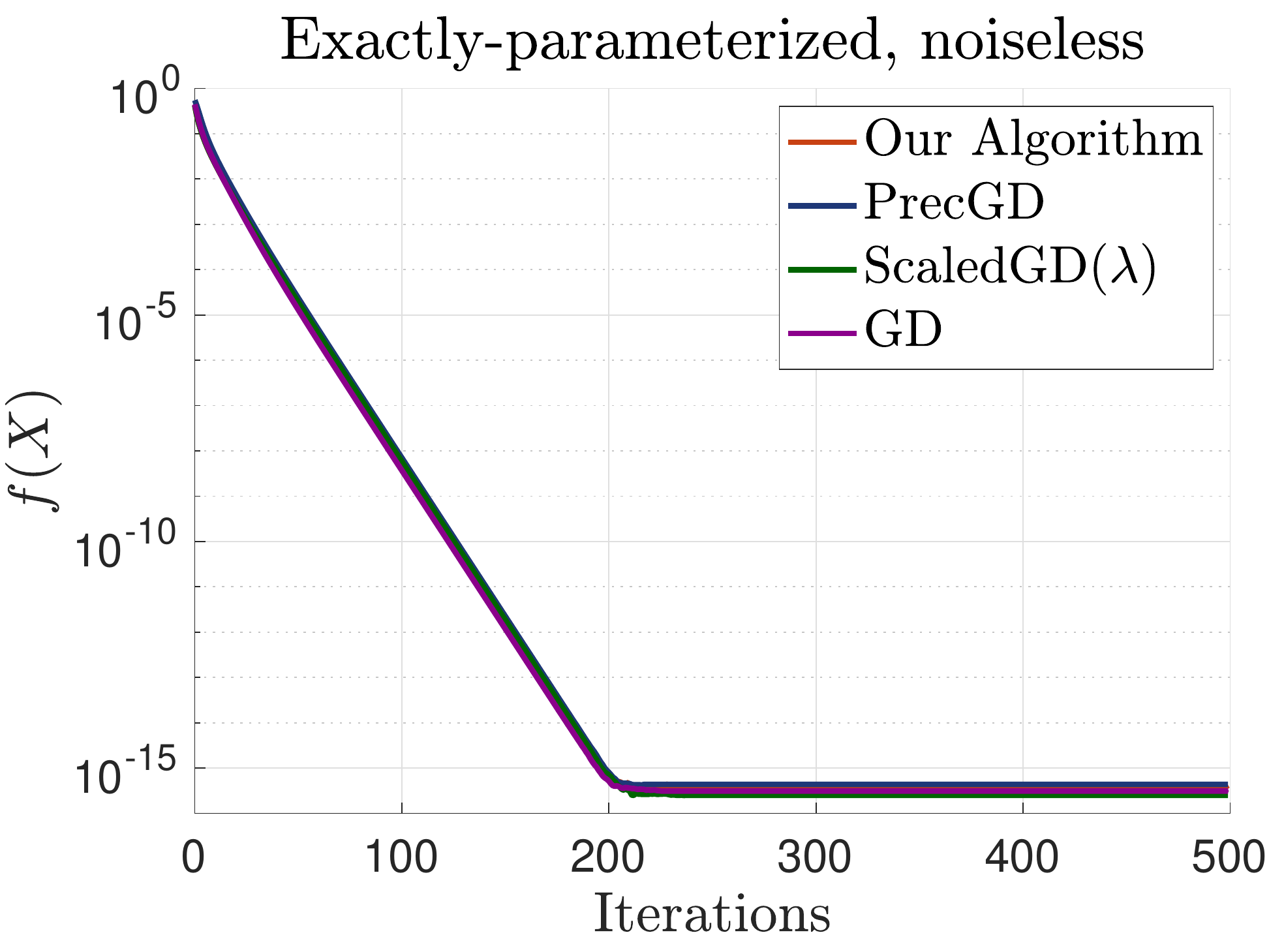}%
	\includegraphics[width = 0.5\textwidth]{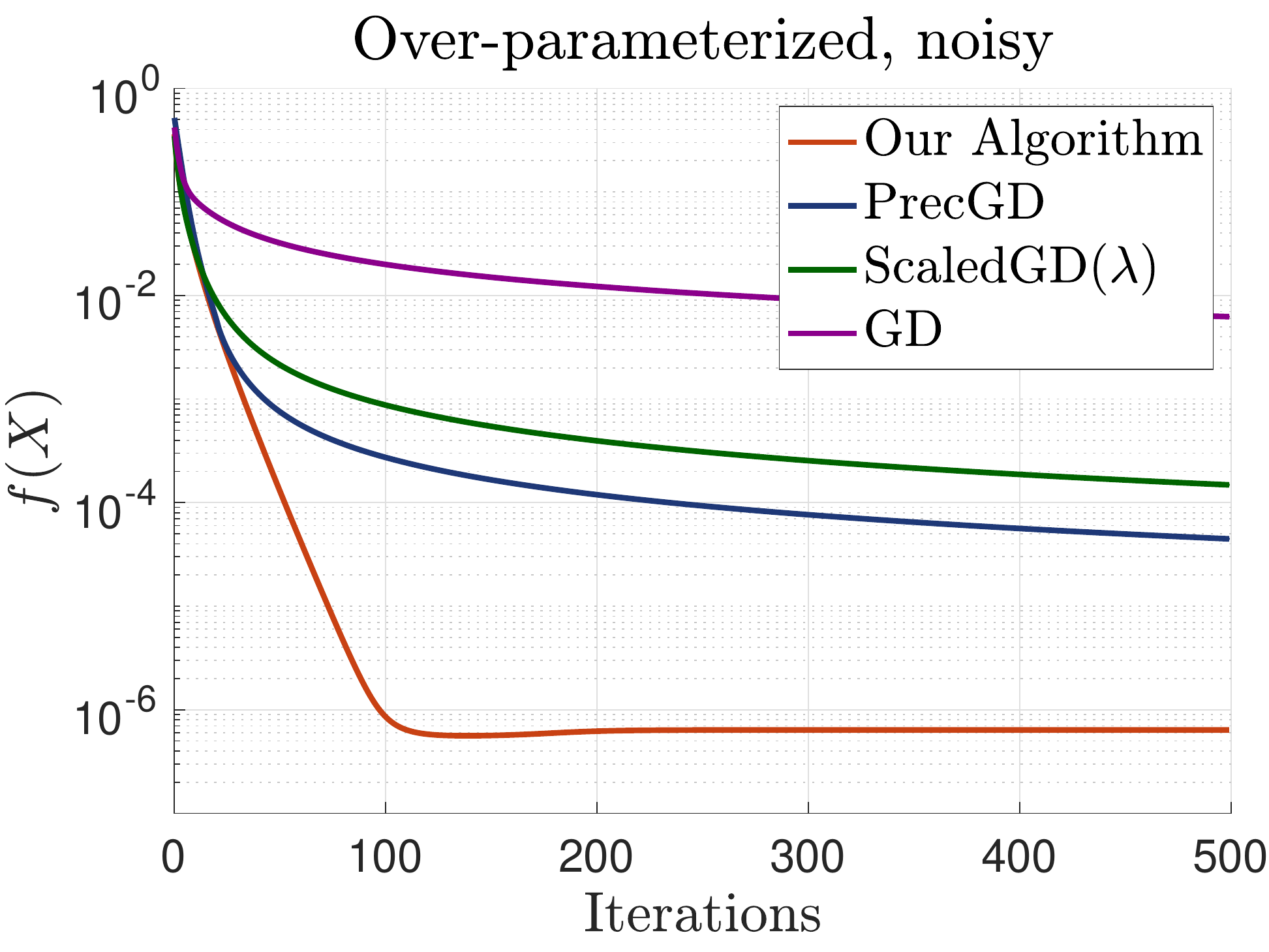}
	\caption{\textbf{Convergence of our algorithm, PrecGD, ScaledGD$(\lambda)$ and GD for Gaussian matrix sensing}. Left: Noiseless measurements with $r=r^\star$. Right: Noisy measurements with $r>r^\star$.}
    \label{fig:gaussian_matrix_sensing}
    \vspace{-0.5em}
\end{figure*}

Figure \ref{fig:gaussian_matrix_sensing} plots the convergence of our algorithm, PrecGD, ScaledGD$(\lambda)$ and GD. The first setting corresponds to the case where $r^\star$ is known, and our measurements are perfect. In this highly unrealistic scenario, we see that the all four methods behave identically, converging linearly to machine error. In the second setting, we see that our algorithm converges to a minimax error of around $10^{-6}$, while PrecGD, ScaledGD$(\lambda)$ and GD struggles to attain the same error. Here the slow down of GD is due to over-parameterization, while the showdown of PrecGD and ScaledGD$(\lambda)$ are due to an inaccurate estimate of scaling parameters. 

\begin{figure*}[h]
	\centering
 \includegraphics[width = 0.5\textwidth]{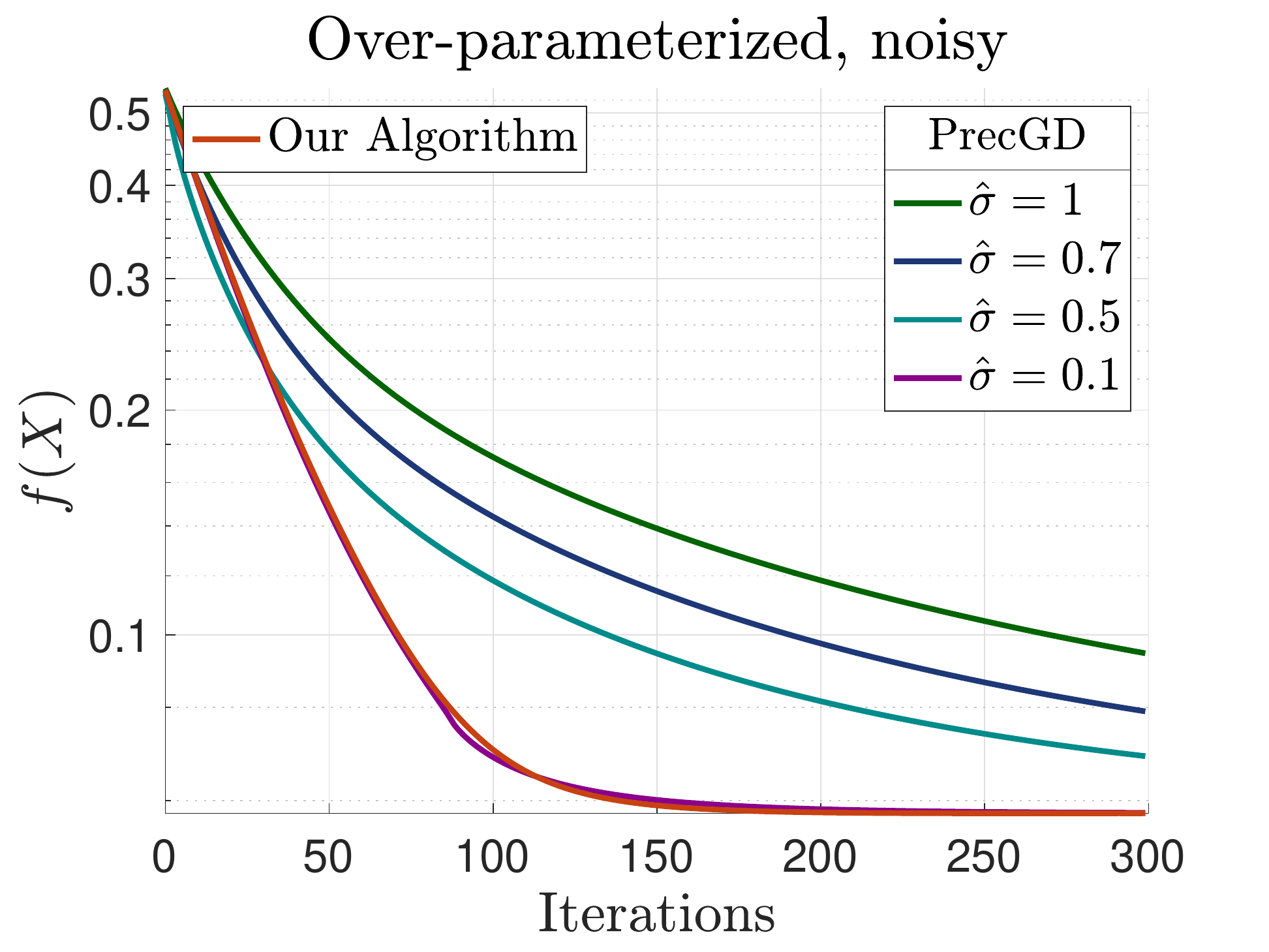}%
	\includegraphics[width = 0.5\textwidth]{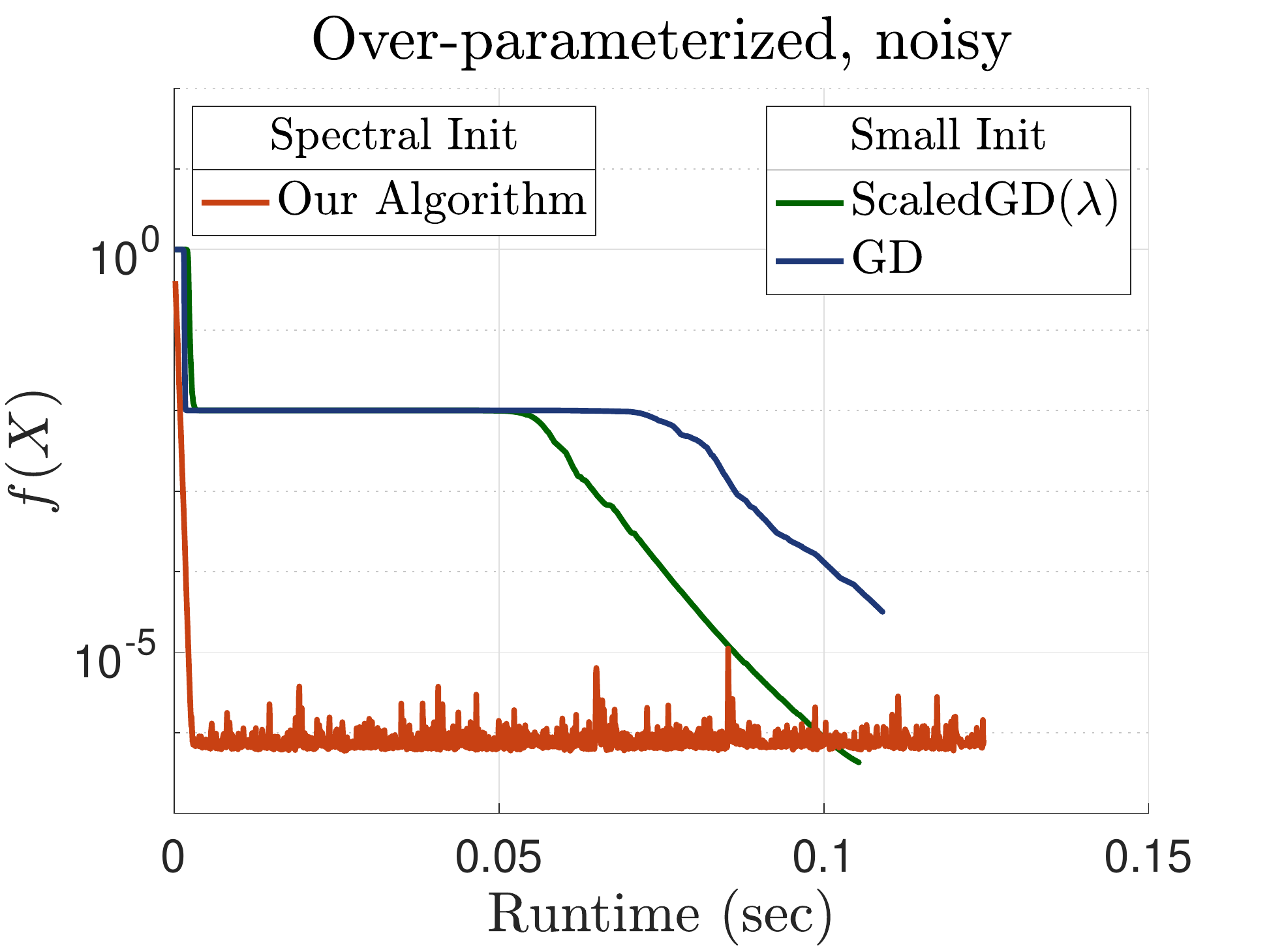}
	\caption{\textbf{Comparison under higher noise setting and small initialization.}. Left: Comparison between our algorithm and PrecGD under higher noise setting. Right: Runtime comparison between our algorithm (with spectral initialization), and ScaledGD$(\lambda)$ and GD (with small initialization).}
    \label{fig:gaussian_matrix_sensing_additional}
    \vspace{-0.5em}
\end{figure*}

\paragraph{High noise setting}
In the first plot of Figure~\ref{fig:gaussian_matrix_sensing_additional}, we plot the convergence of our algorithm under higher noise setting. In particular, we corrupt the measurements with noise $\epsilon_i\sim\mathcal N(0,10^{-1})$. To accommodate higher noise, we set $X$ to be a size $10\times 8$ matrix and minimize $f(X)$ using our algorithm and PrecGD for 500 iterations. In this experiment, our algorithm is implemented with $\beta=0.97$ and PrecGD is implemented with four different proxy variance $\hat \sigma=1,0.7,0.5$ and $0.1$ so that $\eta_t=\sqrt{|f(X_t)-\hat \sigma^2|}$. The learning rate for both methods are set to $\alpha=0.01$. From Figure~\ref{fig:gaussian_matrix_sensing_additional}, we again see that our algorithm converges linearly to the minmax error while PrecGD slows down when the proxy variance is incorrectly estimated.

\paragraph{Comparison with small initialization}
In the second plot of Figure~\ref{fig:gaussian_matrix_sensing_additional}, we compare the runtime of our algorithm (which is initialized using spectral initialization) against ScaledGD$(\lambda)$ and GD that are initialized using small initialization. We note that we include the time for calculating the spectral initial point into the runtime of our algorithm. In this experiment, we set $X$ to be a size $10\times 8$ matrix and minimize $f(X)$ using our algorithm, ScaledGD$(\lambda)$ and GD for around $0.5$ seconds. Here, our algorithm is implemented with $\beta=0.5$, and ScaledGD$(\lambda)$ is implemented with $\lambda=0.01$. The learning rate for all three methods are set to $\alpha=0.1$. In Figure~\ref{fig:gaussian_matrix_sensing_additional}, we see that, despite having to spend extra time to computing the spectral initial points, our algorithm is still significantly faster than ScaledGD$(\lambda)$ and GD.

\subsection{1-bit matrix sensing}
The goal of 1-bit matrix sensing is to recover a ground truth matrix $M^\star$ from 1-bit measurements of each entry in $M^\star$. In particular, the measurements $y_{ij}$ of each entry $M_{ij}^\star$ are quantized, so that they are $1$ with some probability $\sigma(M_{ij}^\star)$ and $0$ with probability $1-\sigma(M_{ij}^\star)$ where $\sigma(\cdot)$ is the sigmoid function. In our experiment on a size $10\times 10$ ground truth matrix $M^\star$, we measure each $y_{ij}$ for a number of times and let $\alpha_{ij}$ denote the percentage of $y_{ij}$ that is equal to 1. To recover the ground truth matrix $M^\star$, we substitute $M^\star = XX^T$ and minimize the following loss function
\[
f(X)= \frac{1}{100}\sum_{i=1}^{10}\sum_{j=1}^{10} -\alpha_{ij}\log\left(\sigma(x_i^Tx_j)\right) - (1-\alpha_{ij})\log\left(1-\sigma(x_i^Tx_j)\right)
\]
where $x_i^T$ is the $i$-th row in $X$, i.e. $X=[x_1\ldots x_{10}]^T$. We perform 1-bit matrix sensing under two settings: the exactly-parameterized, noiseless case; and the over-parameterized, noisy case.

\paragraph{The exactly-parameterized, noiseless case} Recall that the truth rank of $M^\star$ is 2. In the exactly-parameterized case, we set $X$ to be a size $10\times 2$ matrix and minimize $f(X)$ using our algorithm, PrecGD, ScaledGD$(\lambda)$ and GD for 200 iterations. We set $\beta=0.4$ in our algorithm, and $\lambda=0$ in ScaledGD$(\lambda)$. The learning rate for all four methods are set to $\alpha=1$.

\paragraph{The over-parameterized, noisy case} In this setting, we corrupt the measurements with noise $\varepsilon_{ij}\sim\mathcal N(0,10^{-6})$ such that $y_{ij}=1$ with probability $\sigma(M_{ij}^\star+\varepsilon_{ij})$ and $y_{ij}=0$ with probability $1-\sigma(M_{ij}^\star+\varepsilon_{ij})$. We set $X$ to be a size $10\times 4$ matrix and minimize $f(X)$ using our algorithm, PrecGD, ScaledGD$(\lambda)$ and GD for 200 iterations. Here, our algorithm is implemented with $\beta=0.4$, PrecGD is implemented with proxy variance $\hat \sigma=10^{-5}$ so that $\eta_t=\sqrt{|f(X_t)-\hat \sigma^2|}$, and ScaledGD$(\lambda)$ is implemented with $\lambda=0.01$. The learning rate for all four methods are set to $\alpha=1$.

\begin{figure*}[h]
	\centering
	\includegraphics[width = 0.5\textwidth]{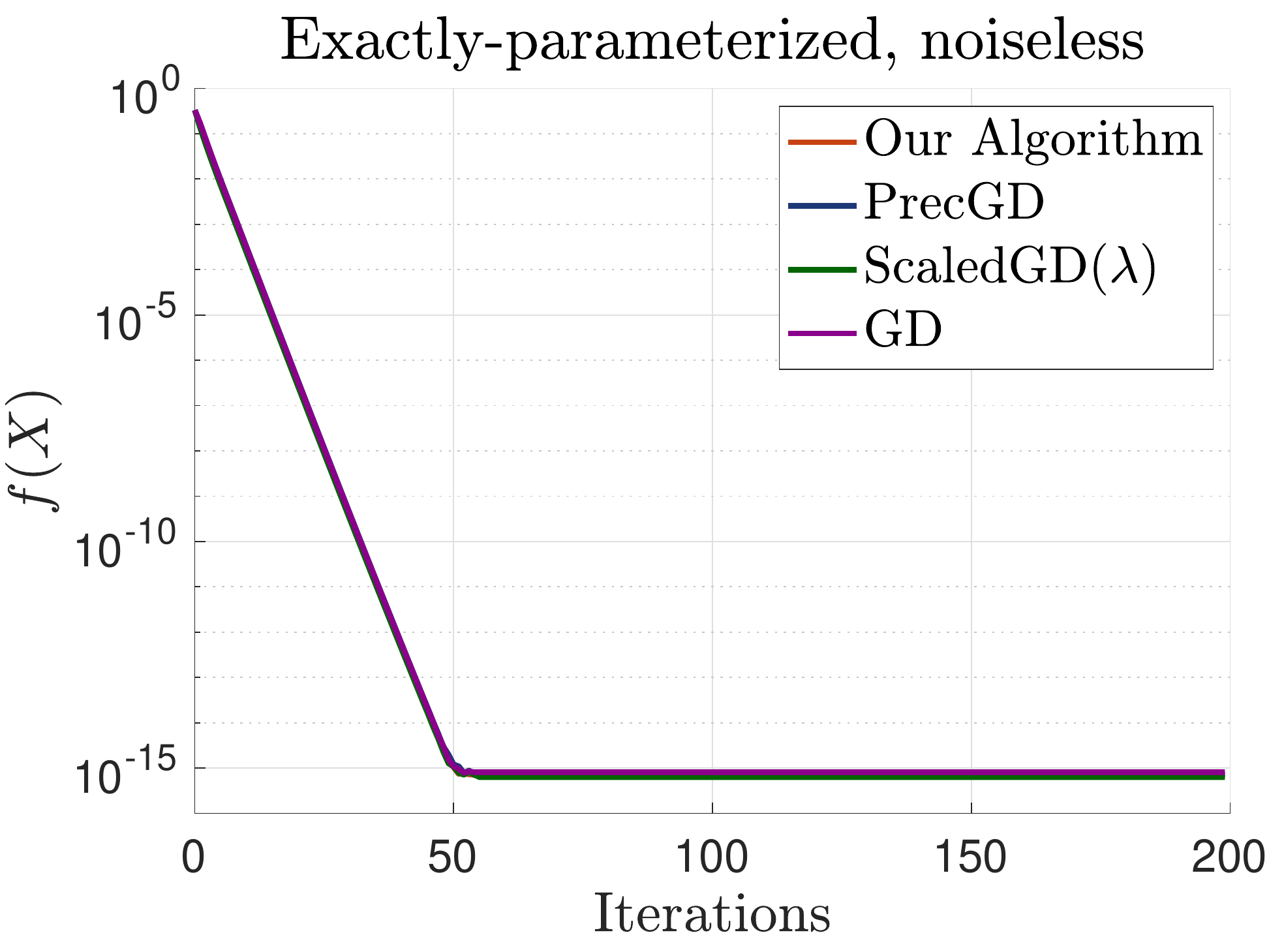}%
	\includegraphics[width = 0.5\textwidth]{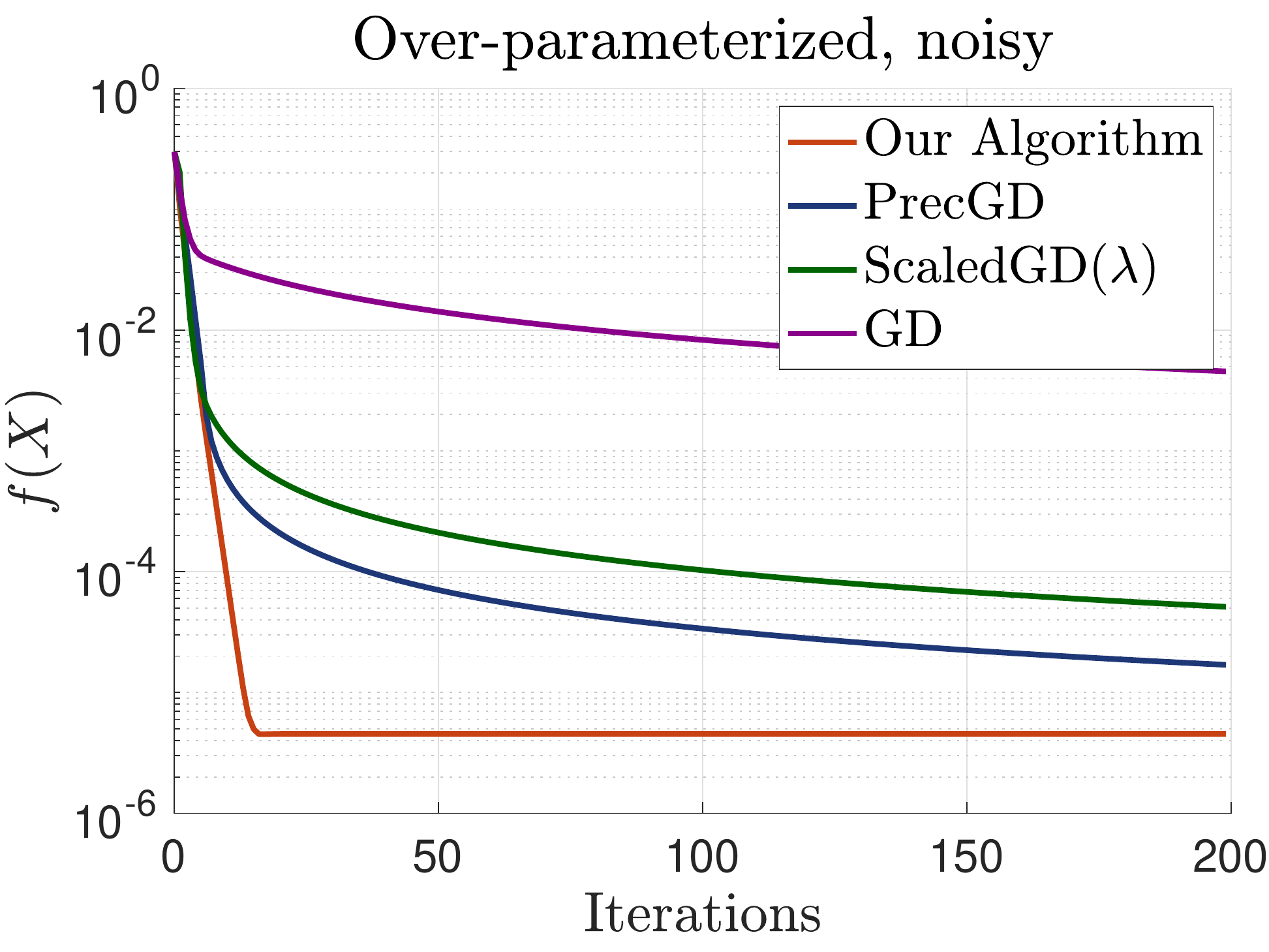}
	\caption{\textbf{Convergence of our algorithm, PrecGD, ScaledGD$(\lambda)$ and GD for 1-bit matrix sensing}. Left: Noiseless measurements with $r=r^\star$. Right: Noisy measurements with $r>r^\star$.}
    \label{fig:1bit_matrix_sensing}
    \vspace{-0.5em}
\end{figure*}

From Figure~\ref{fig:1bit_matrix_sensing}, we again see that the results are almost identical to that of Figure~\ref{fig:gaussian_matrix_sensing}: our algorithm is able to converge linearly to a minimax error rate as soon as $r>r^\star$, but PrecGD, ScaledGD$(\lambda)$ and GD showdown dramatically.

\subsection{Phase retrieval} The goal of phase retrieval is to recover a vector $z\in\mathcal C^n$ from the phaseless measurements of the form $y_i=|\inp{a_i}{z}|^2$ where $a_i\in\mathcal C^n$ are the measurement vectors.  Equivalently, we can view this problem as recovering a complex matrix $M^\star$ from measurements $y_i=\inp{a_ia_i^T}{M^\star}$, subjecting to a constraint that $M^\star$ is rank-1. In our experiment on a length 10 ground truth vector $z$, we set $M^\star=zz^T$ and take 80 measurements on $M^\star\in\mathcal C^{10\times 10}$ using 80 linearly independent measurement vectors $a_i\in\mathcal C^{10}$ drawn from standard Gaussian.  Substituting $M^\star=XX^T$, the loss function of phase retrieval is defined as
\[
f(X)=\frac{1}{80}\sum_{i=1}^{80}\left(\inp{a_ia_i^T}{XX^T}-y_i\right)^2.
\]
We again perform phase retrieval under two settings: the exactly-parameterized, noiseless case; and the over-parameterized, noisy case.

\paragraph{The exactly-parameterized, noiseless case} Recall that the truth rank of $M^\star$ is 1. In the exactly-parameterized case, we set $X$ to be a size $10\times 1$ complex matrix and minimize $f(X)$ using our algorithm, PrecGD, ScaledGD$(\lambda)$ and GD for 1000 iterations. We set $\beta=0.1$ in our algorithm, and $\lambda=0$ in ScaledGD$(\lambda)$. The learning rate for all four methods are set to $\alpha=0.02$.

\paragraph{The over-parameterized, noisy case} In this setting, we corrupt the measurements with noise $\varepsilon_{i}\sim\mathcal N(0,10^{-6})$ such that $y_{i}=\inp{a_ia_i^T}{M^\star}+\varepsilon_{i}$. We set $X$ to be a size $10\times 2$ matrix and minimize $f(X)$ using our algorithm, PrecGD, ScaledGD$(\lambda)$ and GD for 1000 iterations. Here, our algorithm is implemented with $\beta=0.1$, PrecGD is implemented with proxy variance $\hat \sigma=10^{-5}$ so that $\eta_t=\sqrt{|f(X_t)-\hat \sigma^2|}$, and ScaledGD$(\lambda)$ is implemented with $\lambda=0.01$. The learning rate for all four methods are set to $\alpha=0.02$. 

\begin{figure*}[h]
	\centering
	\includegraphics[width = 0.5\textwidth]{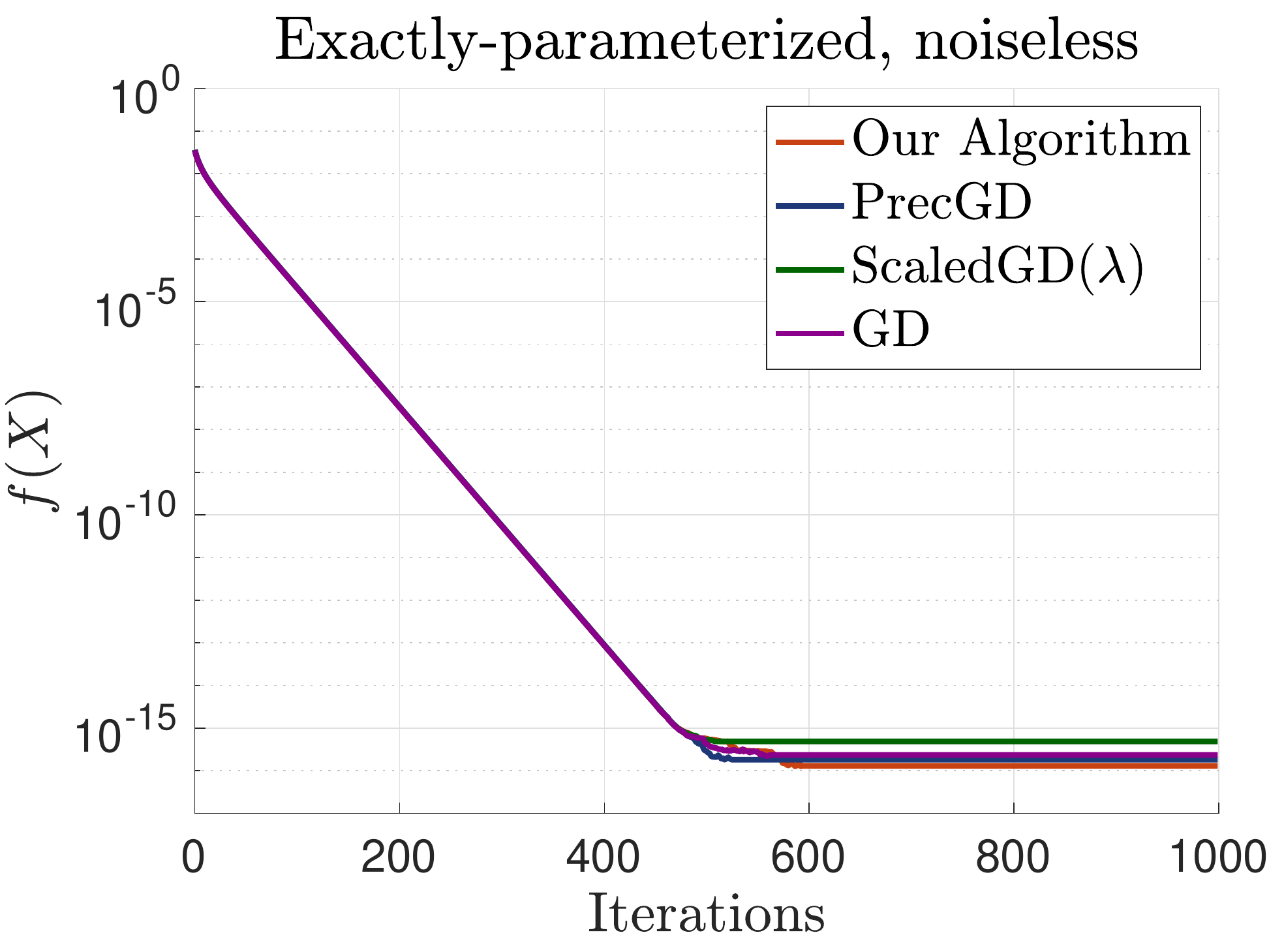}%
	\includegraphics[width = 0.5\textwidth]{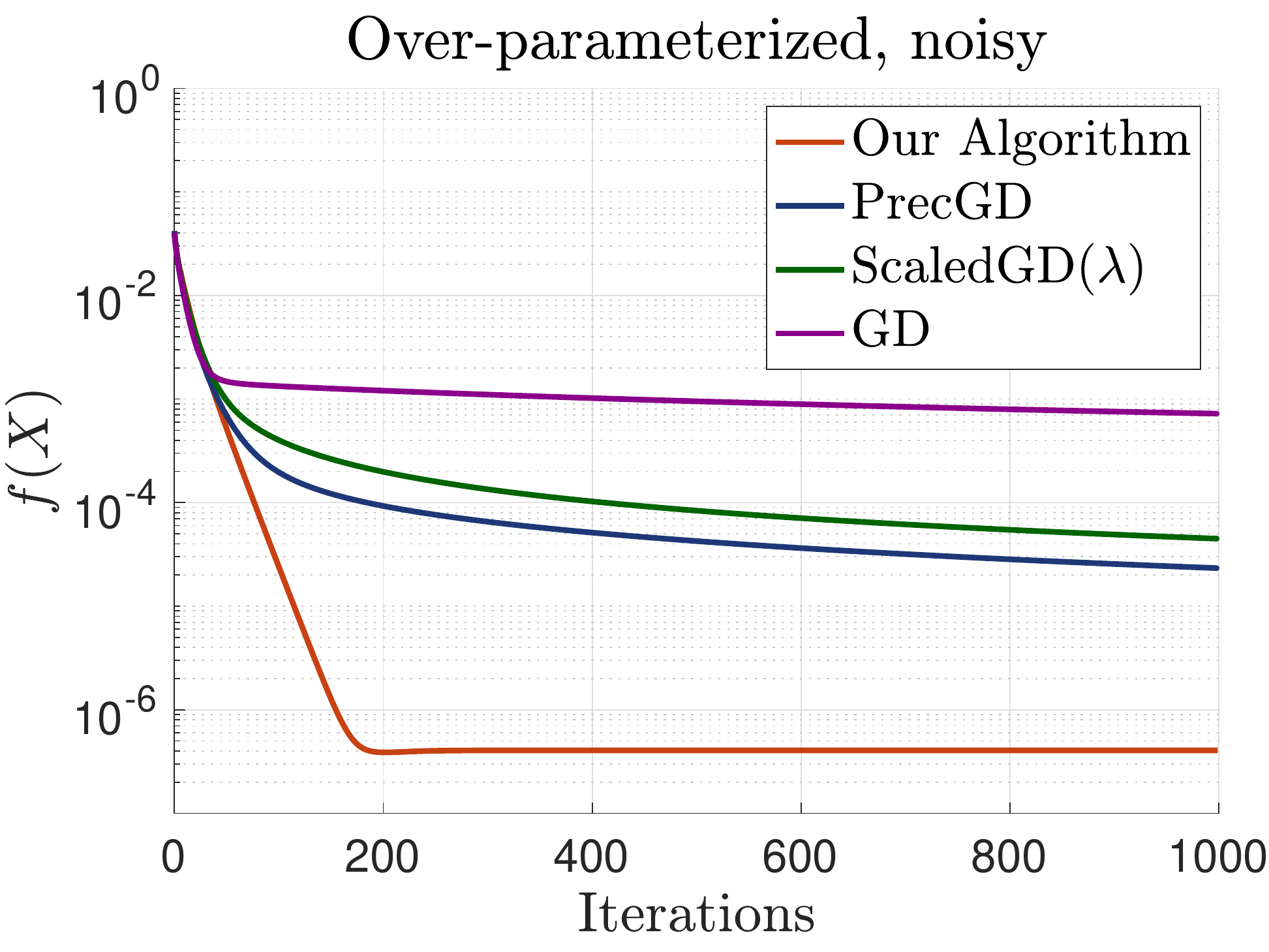}
	\caption{\textbf{Convergence of our algorithm, PrecGD, ScaledGD$(\lambda)$ and GD for phase retrieval}. Left: Noiseless measurements with $r=r^\star$. Right: Noisy measurements with $r>r^\star$. }
    \label{fig:phase_retrieval}
    \vspace{-0.5em}
\end{figure*}

Figure~\ref{fig:phase_retrieval} shows the convergence of our alogirthm, PrecGD, ScaledGD$(\lambda)$ and GD. Again, our algorithm again converge linearly to the minimax error. 

\section{Additional Experiments on Ultrasound Image Recovery}
\begin{figure}
	\centering
    \includegraphics[width = 0.9\textwidth]{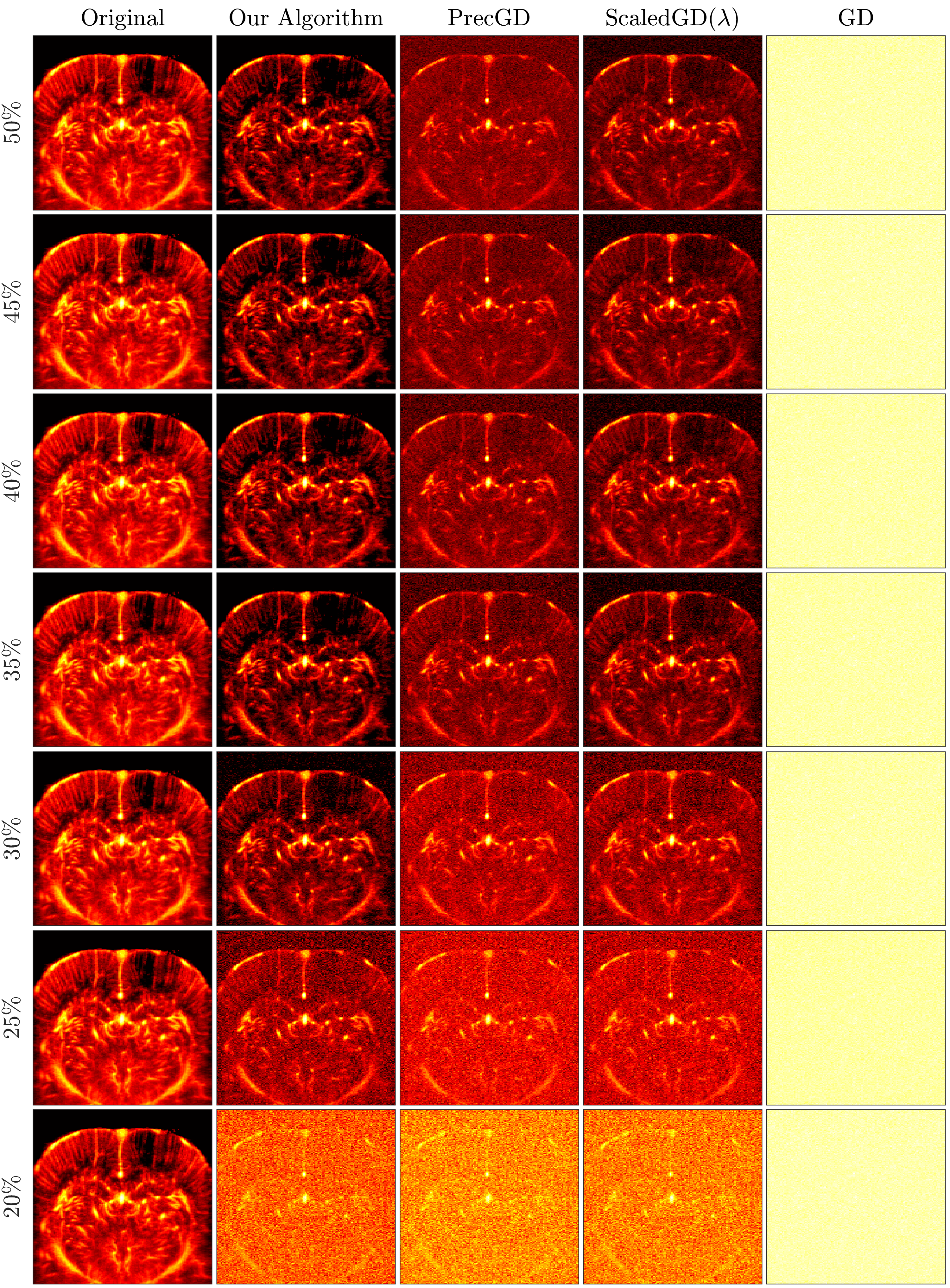}
	\caption{\small\textbf{Denoising an ultrafast ultrasound image under different downsampling rate}. We denoise the ultrasound image in Figure~\ref{fig:ultrasound} under 7 downsampling rates: 50\%, 45\%, 40\%, 35\%, 30\%, 25\% and 20\%. The ultrasound images are shown using power Doppler \cite{bercoff2011ultrafast}. Column 1: original image. Column 2: image denoised from our algorithm \eqref{alg:1}. Column 3: image denoised from ScaledGD$(\lambda)$. Column 4: image denoised from PrecGD. Column 5: image denoised from GD.} 
    \label{fig:ultrasound_sampling}
\end{figure}

We repeat the experiment on ultrasound image denoising task under 7 downsampling rates: 50\%, 45\%, 40\%, 35\%, 30\%, 25\% and 20\%. In all 7 cases, we set the search rank to be $r=100$ and apply our algorithm, PrecGD, ScaledGD$(\lambda)$ and GD to minimize the corresponding loss function $f(U,V)$ for 30 iterations. Our algorithm, PrecGD and GD are implemented using the same hyperparameters and learning rates in Figure~\ref{fig:ultrasound}. ScaledGD$(\lambda)$ is implemented with $\lambda=5\times 10^{-2}$ and learning rate $\alpha=10^7$. As shown in Figure~\ref{fig:ultrasound_sampling}, our algorithm is able to almost perfectly denoise the ultrasound image when the downsampling rate is above 25\%.

\end{document}